\documentclass{amsart}

\usepackage{amssymb}
\usepackage[dvips]{graphics}

\newcounter{internalclaimcounter}

\newtheorem{theorem}{Theorem}[section]
\newtheorem{corollary}[theorem]{Corollary}
\newtheorem{proposition}[theorem]{Proposition}
\newtheorem{lemma}[theorem]{Lemma}
\newtheorem{internalclaim}[internalclaimcounter]{Claim}

\newtheorem*{claim*}{Claim}
\newtheorem*{subclaim}{Subclaim}

\theoremstyle{definition}
\newtheorem{definition}[theorem]{Definition}

\theoremstyle{remark}
\newtheorem{remark}{Remark}[section]

\numberwithin{equation}{section}

\newcommand{\Diagram}{\mathcal{D}}
\newcommand{\Strip}{\mathcal{S}}
\newcommand{\Twist}{\tau}

\newcommand{\Remainder}{\rho}
\newcommand{\Cancel}{\kappa}
\newcommand{\Diam}{\mathrm{Diam}}
\newcommand{\Edge}{e}
\newcommand{\Edgepath}{\gamma}
\newcommand{\EdgepathSystem}{\Gamma}
\newcommand{\BasicEdgepath}{\lambda}
\newcommand{\BasicEdgepathSystem}{\Lambda}
\newcommand{\ConstantEdgepaths}{\Gamma_{\mathrm{const}}}

\newcommand{\const}{\mathrm{const}}
\newcommand{\dec}{\mathrm{dec}}
\newcommand{\inc}{\mathrm{inc}}
\newcommand{\Crossing}{\mathrm{cr}}

\newcommand{\Tangle}{T}
\newcommand{\NumTangles}{N}

\newcommand{\Slope}{R}

\newcommand{\typeI}{\mathrm{I}}
\newcommand{\typeII}{\mathrm{II}}
\newcommand{\typeIII}{\mathrm{III}}
\newcommand{\classA}{\mathrm{A}}
\newcommand{\classB}{\mathrm{B}}

\newcommand{\angleb}[1]{\langle #1 \rangle} 

\newcommand{\circleb}[1]{\langle #1 \rangle^{\circ}} 

\newcommand{\func}{\varphi}

\newcommand{\notionNotInHO}{$ ^\dag$}

\newcommand{\parestart}{(}
\newcommand{\pareend}{)}

\begin{document}

\title{Lower bounds on boundary slope diameters for Montesinos knots}

\author{Kazuhiro Ichihara 
}
\address{%
College of General Education, 
Osaka Sangyo University, 
3--1--1 Nakagaito, Daito, Osaka 574--8530, Japan }
\email{ichihara@las.osaka-sandai.ac.jp}
\thanks{%
  The first author is partially supported by
  Grant-in-Aid for Young Scientists (B), No. 18740038,
  Ministry of Education, Culture, Sports, Science and Technology, Japan.
}

\author{%
    Shigeru Mizushima}
\address{%
        Department of Mathematical and Computing Sciences, 
        Tokyo Institute of Technology, 
        2--12--1 Ohokayama, Meguro, 
        Tokyo 152--8552, Japan}
\email{mizusima@is.titech.ac.jp}

\keywords{boundary slopes, diameter, Montesinos knot}
\subjclass[2000]{Primary 57M25}

\date{\today}


\begin{abstract}
In this paper, two lower bounds on the diameters of the boundary slope sets 
are given for Montesinos knots. 
One is described in terms of the minimal crossing numbers of the knots, 
and the other is related to the Euler characteristics of essential surfaces 
with the maximal/minimal boundary slopes. 
\end{abstract}

\maketitle


\section{Introduction}


Consider a compact possibly non-orientable surface
properly embedded in a knot exterior in the 3-sphere $S^3$. 
It is called \textit{essential} if it is incompressible and boundary-incompressible. 
The boundary of an essential surface 
consists of a parallel family of non-trivial simple closed curves. 
Thus,  on the peripheral torus of the knot, they determine an isotopy class of 
a non-trivial unoriented simple closed curve, 
which is called the \textit{boundary slope} of the surface.


In \cite{H}, Hatcher showed that such boundary slopes are only finitely many. 
Moreover, by Culler and Shalen \cite{CS84}, 
it was proved that there always exist at least two such boundary slopes. 
Therefore one can get a non-empty, finite set of boundary slopes for a knot, 
which is said to be the \textit{boundary slope set}. 

The boundary slope set for a knot $K$ in $S^3$ 
gives a non-empty, finite subset of rational numbers in a standard way. 
See \cite{Ro} for example. 
In view of this, Culler and Shalen introduced and studied in \cite{CS99} 
the \textit{diameter} of the boundary slope set for $K$, 
which is defined as the difference 
between the maximum and the minimum among non-meridional elements. 
This will be denoted by $\Diam(K)$ in this paper.

\subsection{}

Actually Culler and Shalen proved in \cite{CS99} that; 
\[
\Diam(K) \ge 2
\]
holds for a non-trivial knot $K$ in $S^3$ 
if $K$ does not have meridional boundary slope.


This can be generalized for an alternating knot $K$ 
by using the result in \cite{A, DR} as follows: 
\[
\Diam(K) \ge 2\,\Crossing(K)
\]
where $\Crossing(K)$ denotes the minimal crossing number of $K$. 
See \cite{IM2}.


In this paper, we consider a \textit{Montesinos knot}, 
which is obtained by composing a number of rational tangles in line, 
and give the following: 
%
\begin{theorem}
\label{Thm:Diam:LowerBound:ByCrossing}
For a Montesinos knot $K$, we have the inequality 
\[
\Diam(K) \ge 2\,\Crossing(K) -6 \, .
\]
\end{theorem}


This theorem together with the results in \cite{IM2} and 
\cite{MMR} for two-bridge knots 
gives the following upper and lower bounds. 
%
%
\begin{corollary}
\label{Cor:Diam:LowerBound:ByCrossing}
For a non-trivial Montesinos knot K,
we have the inequalities 
\[
2\,\Crossing(K) -6 \le \Diam(K) \le 2\,\Crossing(K) \, . 
\]
\end{corollary}

In fact, as claimed in \cite{IM2}, 
if $K$ is an alternating Montesinos knot, we have the equality 
\[
\Diam(K)=2\,\Crossing(K) \, .
\]

\subsection{}


For an alternating knot $K$, 
the lower bound of $ \Diam(K)$ in terms of $\Crossing (K)$ 
is achieved by considering the checkerboard surfaces $F_1 , F_2$ 
for its reduced alternating diagram. 
In fact, these are known to be essential by \cite{A, DR}. 
For these surfaces, the following inequality also holds; 
\[
\Diam (K) \ge 2\,( (-\chi (F_1) )+ ( -\chi (F_2) ) + 4 
\]

We also generalize this for Montesinos knots as follows:

\begin{theorem}
\label{Thm:Diam:LowerBound:Euler}
Let $K$ be a Montesinos knot. 
Among its non-meridional boundary slopes, 
let $\Slope_1$ and $\Slope_2$ be the maximum and the minimum respectively. 
Then 
there exist two essential surfaces $F_1$ and $F_2$ 
with boundary slopes $\Slope_1$ and $\Slope_2$ 
such that 
\begin{eqnarray}
\Diam (K) = |\Slope_1-\Slope_2|&\ge& 2\,\left( \frac{-\chi}{\sharp s} (F_1) + \frac{-\chi}{\sharp s} (F_2) \right) \, ,
\label{Eq:Diam:LowerBound:Main}
\end{eqnarray}
where 
$ \frac{- \chi}{\sharp s}  (F_i) $ denotes 
the ratio of the negative of the Euler characteristic 
and the number of sheets for $F_i$ for $i=1,2$. 
\end{theorem}
Here, following \cite{HT}, by the \textit{number of sheets} of an essential surface, 
we mean 
the minimal number of intersection between the surface and the meridian of the knot.


\begin{remark}		
We remark that our lower bound is optimal in a sense. 
See Remark \ref{Rem:Optimal}. 
Also it should be compared with 
an upper bound given in \cite[Theorem 4]{IM1}: 
We actually showed that 
\begin{eqnarray*}
|\Slope_1-\Slope_2|&\le& 2\,\left( \frac{-\chi}{\sharp s} (F_1)+\frac{-\chi}{\sharp s} (F_2) \right) + 4 .
\end{eqnarray*}
holds for any pair of essential surfaces $F_1, F_2$ with boundary slopes $\Slope_1, \Slope_2$ in a Montesinos knot exterior. 
\end{remark}


The minimal geometric intersection number of 
the curves representing two slopes $R_1$ and $R_2$ 
is called the \textit{distance} of $R_1$ and $R_2$. 
This is usually denoted by $\Delta ( R_1  , R_2 )$. 
About the distance of the maximal and minimal boundary slopes for a Montesinos knot, 
we have the following 
as an immediate corollary to Theorem \ref{Thm:Diam:LowerBound:Euler}. 
%
\begin{corollary}
\label{Cor:Diam:LowerBound:OnDistance:ByEulerChara}
Let $K$ be a non-trivial Montesinos knot. 
Among its non-meridional boundary slopes, 
let $\Slope_1$ and $\Slope_2$ be the maximum and the minimum respectively. 
Then 
there exist two essential surfaces $F_1$ and $F_2$ 
with boundary slopes $\Slope_1$ and $\Slope_2$ 
such that 
\begin{eqnarray}
\Delta ( \Slope_1 , \Slope_2 ) &\ge& 2\,\left( \frac{-\chi}{\sharp b} (F_1) + \frac{-\chi}{\sharp b} (F_2) \right) \, ,
\end{eqnarray}
where 
$ \frac{- \chi}{\sharp b}  (F_i) $ denote 
the ratio of the negative of the Euler characteristic 
and the number of boundary components for $F_i$ for $i=1,2$. 
\end{corollary}
\begin{proof}
Recall that the distance of the slopes $p/q$ and $r/s$ is calculated by $| ps-qr|$. 
Also note that, for an essential surface $F$ with boundary slope $R$, 
the number of sheets of $F$ is equal to 
the product of the denominator of $R$ with 
the number of boundary components of $F$. 
With these facts, 
the corollary follows from Theorem \ref{Thm:Diam:LowerBound:Euler} immediately. 
\end{proof}
This can be regarded as 
a generalization of the following result shown by Culler and Shalen: 
Suppose that $M$ is a non-exceptional two-surface knot manifold: 
That is, $M$ is 
an irreducible, connected, compact, orientable 3 -manifold with single torus boundary 
such that 
it has at most two distinct isotopy classes of strict essential surfaces and 
it is neither Seifert fibered nor an exceptional graph manifold. 
Let $F_1$ and $F_2$ be 
representatives of the two isotopy classes of connected strict essential surfaces. 
Let $R_i$ denote the boundary slope of $F_i$ and 
let $\sharp b_i$ denote the number of boundary components of $F_i$. 
Then for $i = 1, 2$ we have 
$$
\Delta ( R_1 , R_2 ) \geq 2 \frac{ - \chi (F_i )}{ \sharp b_1 \cdot \sharp b_2 } \ .
$$ 
Please see \cite{CS04} for details.


This paper is organized as follows. 
In the next section, we review the algorithm by Hatcher and Oertel given in \cite{HO}.
In Section \ref{Sec:Twist},
we give some formulae to calculate the twist
and prove Theorem \ref{Thm:Diam:LowerBound:ByCrossing}.
In the last section,
we introduce the remainder term,
give formulae for calculating its value
and 
show Theorem \ref{Thm:Diam:LowerBound:Euler}.

\section{Montesinos knots and Algorithm of Hatcher-Oertel}
\label{Sec:HO}

In this section, 
we give a brief review of the Hatcher-Oertel's work given in \cite{HO}, 
which is the base of our arguments. 
We also prepare basic terminologies used in the rest of the paper. 
Note that the terms marked with ``\dag'' are about notions 
which do not appear in \cite{HO}
and are introduced by the authors in accordance with our argument.

\subsection{Montesinos knot}
\label{Subsec:HO:Montesinos knot}

\begin{figure}[htb]
 \begin{center}
  \begin{picture}(109,91)
   \put(0,0){\scalebox{0.25}{\includegraphics{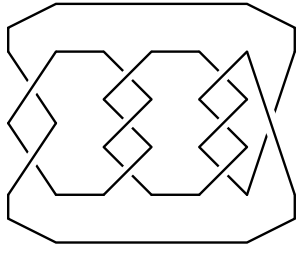}}}
  \end{picture}
  \caption{A diagram of $M(1/2,1/3,-2/3)$}
  \label{Fig:MontesinosKnot}
 \end{center}
\end{figure}

Let us start with the definition of Montesinos knots. 
A \textit{Montesinos knot} is defined as a knot obtained by 
putting rational tangles together in a row. 
See Figure \ref{Fig:MontesinosKnot} for example. 
A Montesinos knot obtained from rational tangles 
$\Tangle_1,\Tangle_2,\ldots,\Tangle_\NumTangles$ 
will be denoted by $M(\Tangle_1,\Tangle_2,\ldots,\Tangle_\NumTangles)$. 
Here and in the sequel, $\Tangle_i$ denotes 
an irreducible fraction or the corresponding rational tangle 
depending on the situation.
In the following, 
we assume that each rational tangle is non-integral, 
just for normalization.
Furthermore, we will always assume that the number of tangles is at least three. 
Note that the knots with at most two tangles are two-bridge knots. 
For two-bridge knots, 
Theorem \ref{Thm:Diam:LowerBound:ByCrossing} holds by \cite{MMR}. 
Also Theorem \ref{Thm:Diam:LowerBound:Euler} holds since 
non-trivial two-bridge knots are all alternating. 

\subsection{Hatcher-Oertel's algorithm}

We here give a very brief review of the algorithm of Hatcher and Oertel. 
See \cite{HO} or \cite{IM1} for detail.

Given Montesinos knot $K=M(\Tangle_1,\Tangle_2,\ldots,\Tangle_\NumTangles)$, 
divide the 3-sphere $S^3$ into $\NumTangles$ 3-balls such that 
$K$ is decomposed into $\NumTangles$ rational tangles $\parestart \Tangle_1, \Tangle_2, \ldots, \Tangle_\NumTangles\pareend$. 
At the same time, 
an essential surface $F$ embedded in the exterior of $K$ 
is divided into surfaces $\parestart F_1, F_2, \ldots, F_\NumTangles\pareend$.
Each of these surfaces $F_i$ 
can be isotoped into some standard position,
and is represented by an ``edgepath'' $\Edgepath_i$ in a ``diagram'' $\Diagram$,
and then, 
whole $F$ is represented by an ``edgepath system'' $\EdgepathSystem$, 
which consists of $\NumTangles$ ``edgepaths'' as follows.

\begin{remark}\label{remark:ambiguity}
An ``edge" corresponds to a piece of surface called a ``saddle".
An ``edgepath" consists of edges, 
and so, it corresponds to a surface obtained by combining saddles. 
Now, in fact, there are two possible choices of making saddles for single edge.
Thus, multiple surfaces correspond to an ``edgepath", and also to an ``edgepath system''.
We here remark that, despite this ambiguity,
all the surfaces corresponding to the same ``edgepath system" 
have the common boundary slope,
the common value of $-\chi/\sharp s$, 
the common ``twist", 
and the common ``remainder term'' defined in a later section. 
See \cite{IM1} for details. 
\end{remark}

\subsubsection{Diagram} 
\begin{figure}[htb]

   \begin{minipage}{180pt}
   \begin{center}
    \scalebox{0.75}{
    \begin{picture}(70,150)
     \put(0,0){\scalebox{0.7}{\includegraphics{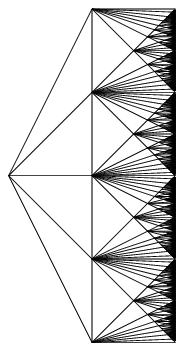}}}
     \put(49,135){\rotatebox{90}{\scalebox{1.0}{$\cdots$}}}
     \put(49,-10){\rotatebox{90}{\scalebox{1.0}{$\cdots$}}}
     \put(2,70.5){\scalebox{1.5}{\vector(1,0){55}}}
     \put(42,5){\scalebox{1.5}{\vector(0,1){90}}}
     \put(29,15){\tiny $-2$}
     \put(29,40){\tiny $-1$}
     \put(36,98){\tiny $1$}
     \put(36,123){\tiny $2$}
     \put(40,148){$v$}
     \put(3,82){\vector(1,-1){10}}
     \put(-10,85){$\angleb{1/0}$}
     \put(2,64){\tiny $-1$}
     \put(72,64){\tiny $1$}
     \put(92,68){$u$}
     \put(33,63){\tiny $O$}
    \end{picture}
    }
   \end{center}
   \caption{The diagram $\Diagram$}
   \label{Fig:Diagram}
   \end{minipage}
\hspace{-10mm}
    \begin{minipage}{200pt}
    \begin{center}    
    \scalebox{0.75}{
     \begin{picture}(150,150)
      \put(0,-10){\scalebox{0.75}{\includegraphics{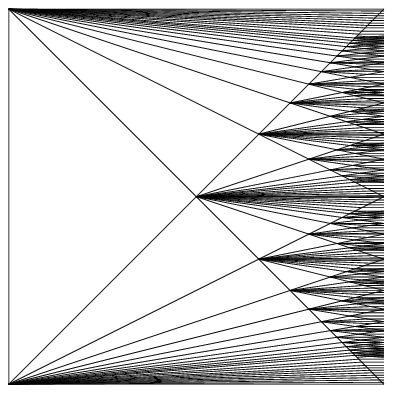}}}
      \put(-5,-2){\vector(1,0){10}}
      \put(-21,-5){$\angleb{0}$}
      \put(156,-2){\vector(-1,0){10}}
      \put(158,-5){$\circleb{0}$}
      \put(-5,133){\vector(1,0){10}}
      \put(-21,130){$\angleb{1}$}
      \put(156,133){\vector(-1,0){10}}
      \put(158,130){$\circleb{1}$}
      \put(63,66){\vector(1,0){10}}
      \put(46,63){$\angleb{\frac{1}{2}}$}
      \put(156,66){\vector(-1,0){10}}
      \put(158,63){$\circleb{\frac{1}{2}}$}
      \put(83,43){\vector(1,0){10}}
      \put(66,40){$\angleb{\frac{1}{3}}$}
      \put(156,43){\vector(-1,0){10}}
      \put(158,40){$\circleb{\frac{1}{3}}$}
      \put(83,88){\vector(1,0){10}}
      \put(66,85){$\angleb{\frac{2}{3}}$}
      \put(156,88){\vector(-1,0){10}}
      \put(158,85){$\circleb{\frac{2}{3}}$}
     \end{picture}
    }
    \end{center}
    \caption{
     A part of the diagram $\Diagram$ in 
     $[0,1]\times[0,1]$
    }
    \label{Fig:Diagram2}
    \end{minipage}

  \end{figure}

The graph on the $u$-$v$ plane defined as follows 
is called a \textit{diagram} and denoted by $\Diagram$. 
A vertex of $\Diagram$ is a point $(u,v)=((q-1)/q,p/q)$ denoted by $\angleb{p/q}$
or a point $(u,v)=(1,p/q)$ denoted by $\circleb{p/q}$
for an irreducible fraction $p/q$ with $q> 0$, 
or a point $(u,v)=(-1,0)$ denoted by $\angleb{1/0}$. 
Two vertices $\angleb{p/q}$ and $\angleb{r/s}$ 
are connected by an edge, which is a straight segment, if $|ps-qr|=1$.
There is another kind of edge called a \textit{horizontal edge},
which connects $\angleb{p/q}$ and $\circleb{p/q}$ for $p/q\neq 1/0$.
An important class of non-horizontal edges is the \textit{vertical edges}, 
which connect the vertices $\angleb{z}$ and $\angleb{z+1}$ for an arbitrary integer $z$.
Another important class of the edges is the \textit{$\infty$-edges},
which connect $\angleb{1/0}$ and $\angleb{z}$ for an integer $z$.
The diagram $\Diagram$ is illustrated 
in Figures \ref{Fig:Diagram} and \ref{Fig:Diagram2}.
%
By $\Strip$, we denote a subgraph of $\Diagram$ lying in $0\le u\le 1$.
%

%
%

\subsubsection{Edgepath and Edgepath system}
A path $\Edgepath$ on $\Diagram$ is called \textit{minimal} 
if $\Edgepath$ is forbidden from 
retracing an edge and going along two sides of single triangle in succession. 
An \textit{edgepath} $\Edgepath$ for a rational tangle 
associated to the irreducible fraction $p/q$ is 
defined as a minimal path on $\Diagram$ such as:
(a) just one point on a horizontal edge $\angleb{p/q}$\,--\,$\circleb{p/q}$, 
or
(b) a path starting at $\angleb{p/q}$ 
proceeding monotonically from right to left or at least vertically. 
An edgepath is called a \textit{constant edgepath} in the former case, 
and a \textit{non-constant edgepath} otherwise. 

An \textit{edgepath system} $\EdgepathSystem$ 
for a Montesinos knot $K=M(T_1,T_2,\ldots,T_\NumTangles)$
is then defined as a tuple 
$\parestart \Edgepath_1,\Edgepath_2,\ldots,\Edgepath_\NumTangles\pareend$ 
of edgepaths 
where each $\Edgepath_i$ is an edgepath for the tangle $T_i$.

\subsubsection{Basic edgepath system and Extended basic edgepath system}
\label{SubsubSec:Basic:ExtendedBasic:EdgepathSystem}
We call a non-constant edgepath for a rational tangle $T_i$ 
whose $u$-coordinate of the endpoint is $0$
and which has no vertical edges a \textit{basic edgepath}\notionNotInHO for $T_i$. 
A \textit{basic edgepath system}\notionNotInHO for a Montesinos knot 
is then defined as a set of basic edgepaths 
with the starting points associated to the tangles. 

A basic edgepath $\BasicEdgepath_i$ for a tangle $T_i$ 
is extended by combining $\BasicEdgepath_i$ and 
the horizontal edge $\angleb{T_i}$\,--\,$\circleb{T_i}$. 
We call it 
the \textit{extended basic edgepath}\notionNotInHO $\widetilde{\BasicEdgepath_i}$ 
of $\BasicEdgepath_i$. 
By an \textit{extended basic edgepath system}\notionNotInHO, 
we mean the set of the extended basic edgepaths 
$\widetilde{\BasicEdgepathSystem}=
\parestart \widetilde{\BasicEdgepath}_1, \widetilde{\BasicEdgepath}_2, \ldots,
\widetilde{\BasicEdgepath}_\NumTangles\pareend$
obtained from a basic edgepath system 
$\BasicEdgepathSystem=\parestart \BasicEdgepath_1, \BasicEdgepath_2, \ldots, \BasicEdgepath_\NumTangles\pareend$.

An extended basic edgepath $\widetilde{\BasicEdgepath}$
is naturally regarded as 
a continuous piecewise-affine function $[0,1]\rightarrow\mathbb{R}$
which gives the $v$ coordinate of the point $(u,v)$ on the edgepath for each $u$. 
We then regard an extended basic edgepath system 
$\widetilde{\BasicEdgepathSystem}$ as a function
defined as the sum 
$\widetilde{\BasicEdgepathSystem}(u)=\sum_{i=1}^{\NumTangles} \widetilde{\BasicEdgepath_i}(u)$
of functions 
$\parestart \widetilde{\BasicEdgepath}_1, \widetilde{\BasicEdgepath}_2, \ldots,
\widetilde{\BasicEdgepath}_\NumTangles\pareend$.
Similarly, a basic edgepath system $\BasicEdgepathSystem$, 
or an edgepath system in general, 
can be regarded as a function, 
which is a restriction of the function $\widetilde{\BasicEdgepathSystem}$.

\subsubsection{Partial edge}
From an extended basic edgepath system 
$\widetilde{\BasicEdgepathSystem}
=
\parestart \widetilde{\BasicEdgepath_1}, \widetilde{\BasicEdgepath_2}, \ldots,
\widetilde{\BasicEdgepath_\NumTangles}\pareend$ 
for a Montesinos knot $K=M(p_1/q_1$, $p_2/q_2$, $\ldots,$ $p_\NumTangles/q_\NumTangles)$
and an arbitrary $u_0$ with $0<u_0<1$, 
by \textit{cutting} it at $u=u_0$, 
we can make an edgepath system 
$\EdgepathSystem=\parestart \Edgepath_1, \Edgepath_2, \ldots, \Edgepath_\NumTangles\pareend$ as follows.
For each $i$,
if $u_0\le(q_i-1)/q_i$, then we set $\Edgepath_i$ to be a non-constant edgepath $\BasicEdgepath_i \cap \{(u,v)|u_0\le u\}=\widetilde{\BasicEdgepath_i} \cap \{(u,v)|u_0\le u\le (q_i-1)/q_i \}$.
If $u_0>(q_i-1)/q_i$, 
then we set $\Edgepath_i$ to be a constant edgepath consisting of 
the point $P_i$ with coordinates $(u,v)=(u_0,p_i/q_i)$ 
on a horizontal edge $\angleb{p_i/q_i}$\,--\,$\circleb{p_i/q_i}$.
In this operation,
the edgepath system $\EdgepathSystem$
may include 
an edge
which is a proper subset of an edge of the diagram $\Diagram$.
The edge is called a \textit{partial edge}.
In contrast,
an edge of the diagram $\Diagram$ completely included in an edgepath
is called a \textit{complete edge}.

\subsubsection{Gluing consistency}
In the way developed in \cite{HO}, 
an edgepath system is associated to a surface 
embedded in a Montesinos knot exterior. 
Such an edgepath system must satisfy the following condition: 
For the endpoints of the edgepaths, 
their $u$-coordinates are common 
and their $v$-coordinates $\parestart v_1, v_2, \ldots, v_\NumTangles\pareend$ are summed up to zero; 
\begin{eqnarray}
\label{Eq:EquationGluingConsistency}
\sum_{i=1}^{\NumTangles} v_i=0
.
\end{eqnarray}
We refer these conditions as the \textit{gluing consistency}, and 
call an edgepath system satisfying these conditions 
a \textit{candidate edgepath system}.

\subsubsection{Type of edgepath system}
In order to enumerate all boundary slopes for a Montesinos knot $K$, 
we first enumerate all basic edgepath systems for $K$, 
which is obviously finitely many, 
and then list candidate edgepath systems up as follows.
(I)
From each basic edgepath system $\BasicEdgepathSystem$, 
by making the extended basic edgepath system $\widetilde{\BasicEdgepathSystem}$, 
solving the equation $\widetilde{\BasicEdgepathSystem}(u)=0$,
getting a solution $u_0>0$, 
and cutting $\widetilde{\BasicEdgepathSystem}$ at $u=u_0$,
we have a candidate edgepath system, 
which is called a \textit{type I} edgepath system. 
The endpoints of its edgepaths have the common $u$-coordinate $u=u_0>0$. 
(II)
By adding appropriate number of vertical edges to some edgepaths of $\BasicEdgepathSystem$ if necessary,
we may have a candidate edgepath system, 
which is called a \textit{type II} edgepath system.
The endpoints of its edgepaths have the common $u$-coordinate $0$. 
(III)
By adding suitable (possibly partial) $\infty$-edges to $\BasicEdgepathSystem$,
we have a candidate edgepath system, 
which is called a \textit{type III} edgepath system.
The endpoints of its edgepaths have the common $u$-coordinate $u<0$. 

Remark here that our classification is slightly different 
from that in \cite{HO}, that is,
a basic edgepath system satisfying gluing consistency 
is regarded as type II in this paper and as type I in \cite{HO}.

By Corollary 2.4 to Proposition 2.10 in \cite{HO},
we can determine whether 
surfaces corresponding to a candidate edgepath system 
are essential or not.
We thus obtain a list of edgepath systems corresponding to some essential surface.
The algorithm of Hatcher and Oertel completes 
by calculating all the boundary slopes for the surfaces 
associated to such edgepath systems.

\subsubsection{Sign of edge}
A non-constant edgepath has a roughly right-to-left direction.
With this direction,
each non-$\infty$ non-horizontal edge in the edgepath
is said to be \textit{increasing} or \textit{decreasing}
according to whether
$v$-coordinate increases or decreases
as a point moves along the edge according to the direction.
For each possibly-partial non-$\infty$-edge $\Edge$ in the edgepath, 
we assign $+1$ or $-1$ as the \textit{sign} $\sigma(\Edge)$
if the edge $\Edge$ is increasing or decreasing.

\subsubsection{Monotonic edgepath system}
At each point $\angleb{p/q}$ with $q\ge 2$,
there exist
exactly one increasing leftward edge and exactly one decreasing leftward edge.
Hence,
for a fixed tangle,
there exists 
the unique basic edgepath 
consisting of only increasing leftward edges,
which is called the \textit{monotonically increasing basic edgepath}\notionNotInHO.
Similarly 
the unique \textit{monotonically decreasing basic edgepath}\notionNotInHO exists.
Note that they are minimal.
Then, for a fixed Montesinos knot $K$,
there exist the unique basic edgepath system
consisting of monotonically increasing basic edgepaths
and 
the unique basic edgepath system
consisting of monotonically decreasing basic edgepaths,
which are called 
the \textit{monotonically increasing basic edgepath system}\notionNotInHO
and 
the \textit{monotonically decreasing basic edgepath system}\notionNotInHO.
Let $\BasicEdgepathSystem_\inc$ and $\BasicEdgepathSystem_\dec$ denote them respectively.
When we regard the corresponding extended basic edgepath systems 
as a function,
naturally,
$\widetilde{\BasicEdgepathSystem_\dec}(u) \le 
\widetilde{\BasicEdgepathSystem}(u) \le
\widetilde{\BasicEdgepathSystem_\inc}(u)$
holds for 
any extended basic edgepath system 
$\widetilde{\BasicEdgepathSystem}$ for a fixed $K$ and any $0\le u\le 1$.
Besides,
$\BasicEdgepathSystem_\inc$ and $\BasicEdgepathSystem_\dec$
are convex and concave as a function respectively.

\subsubsection{Length of edgepath}
\label{SubsubSec:LengthOfEdgepath}
For each possibly partial edge in the edgepath,
we assign the length $|\Edge|$,
where the length of a complete edge is $1$ 
and the length of a partial edge is less than $1$. 
Precisely, as calculated in \cite{HO} or \cite{IM1}, 
a partial edge $e$, 
which is included in a complete edge $\angleb{p/q}$\,--\,$\angleb{r/s}$ 
and has $u_0$ as the $u$-coordinate of the endpoint, has the length 
\begin{eqnarray}
|e|&=& \frac{1+s(u_0-1)}{(s-q)(u_0-1)}
\label{Eq:Formula:LengthOfPartialEdge}
.
\end{eqnarray}

%
\subsection{Incompressibility}
\label{Subsec:HO:Incompressibility}

As mentioned in Remark \ref{remark:ambiguity}, 
a set of surfaces correspond to a candidate edgepath system. 
In view of this, 
a candidate edgepath system is called 
\textit{incompressible}, \textit{compressible} or \textit{indeterminate},
if
all the corresponding surfaces are essential,
all the corresponding surfaces are inessential, 
or
the set of the corresponding surfaces includes both essential and inessential ones, 
respectively.

\begin{remark}\label{remark:pi_1-injectivity}
The aim of the algorithm of Hatcher and Oertel
is to enumerate all the boundary slopes
of the orientable essential surfaces.
In fact, they determine $\pi_1$-injectivity instead of incompressibility.
For instance, the difference of these two notions is mentioned in \cite{HT}.
Note that $\pi_1$-injectivity is stronger than incompressibility.
This does not matter in \cite{HO},
since $\pi_1$-injectivity is equivalent to incompressibility
as for orientable surfaces.
Though, 
since we here deal with possibly non-orientable surfaces,
we have to be careful about the difference.
Nonetheless, in our later argument,
the difference does not make any trouble fortunately.
\end{remark}

Here, we recall two notions used in \cite{HO}, 
and give a lemma about incompressibility of edgepath systems 
used in later arguments.

%
%
Two successive edges 
$\angleb{p_3/q_3}$\,--\,$\angleb{p_2/q_2}$\,--\,$\angleb{p_1/q_1}$ 
on the diagram $\Diagram$ 
are said to be \textit{reversible} 
if the two edges lie in two triangles of $\Diagram$ sharing a common edge. 
For example, $\angleb{1/0}$\,--\,$\angleb{0}$\,--\,$\angleb{1/2}$ is reversible
since these edges lie in two triangles sharing the common edge
$\angleb{0}$\,--\,$\angleb{1}$.
An edgepath is said to be \textit{completely reversible}
if all pairs of successive two edges in it are reversible.

%
%
For an edgepath,
if the last edge of the edgepath is included in $\angleb{p/q}$\,--\,$\angleb{r/s}$,
then the \textit{final $r$-value} of the edgepath is defined to be $(s-q)$.
In some cases, we may give the positive or negative sign to $r$-values
according to whether the last edge is increasing or decreasing respectively.
An edgepath system has a \textit{cycle of final $r$-values} 
obtained by collecting the final $r$-values of $\NumTangles$ edgepaths in the system.

Here, we summarize some results about incompressibility
in propositions in \cite{HO}
used in this paper.
%
\begin{lemma}
 \label{Lem:Incompressibility}
For a Montesinos knot $K=M(\Tangle_1,\Tangle_2,\ldots,\Tangle_\NumTangles)$, 
where $N \ge3$, the following hold. 
 \begin{enumerate}
  \item \label{Lem:Incompressibility:MonotonicTypeI}
        A monotonically decreasing type I edgepath system is incompressible.

  \item \label{Lem:Incompressibility:TypeIorII}
        Suppose that a type I or type II edgepath system $\EdgepathSystem$ lying in $\Strip$ 
        has $(+1,-2,r_3)$ with $r_3\le-5$ as its cycle of final $r$-values. 
	Then $\EdgepathSystem$ is incompressible unless the last edge
	of $\Edgepath_3$ lies in the same triangle of $\Strip$, and has
	the same ending point, as an edge with $r=1$.

  \item \label{Lem:Incompressibility:MonotonicTypeII}
        \textup{(a)} A monotonically decreasing type II edgepath system is incompressible.\\
	\textup{(b)} If a basic edgepath system $\BasicEdgepathSystem$ does not satisfy 
		the condition (*): 
				\begin{quote}
				\upshape
				The cycle of final $r$-values for the $\Edgepath_i$'s contains at least one $-1$,
				and each $-1$ in this cycle is separated from the next $-1$
				by $r_i$'s all but possibly one of which are $-\tilde{2}$'s
			\end{quote}
		then there exists an incompressible or indeterminate type II edgepath system 
		obtained from $\BasicEdgepathSystem$
                by adding upward vertical edges.
		Here 
		an edgepath with $-\tilde{2}$ means
		that the edgepath has final $r$-value $-2$ and is
	completely reversible. 

  \item \label{Lem:Incompressibility:TypeIII}
	A type III edgepath system $\EdgepathSystem$ is compressible
	if and only if 
	$|\EdgepathSystem(0)|\le 1$ holds
	and at least $\NumTangles-2$ edgepaths are completely reversible.
	On the other hand, 
        a type III edgepath system is incompressible 
             if it is constructed 
             from a basic edgepath system $\BasicEdgepathSystem$
             with $|\BasicEdgepathSystem(0)|\ge 2$.
 \end{enumerate}
\end{lemma}
\begin{proof}
(1), (3)(a): 
By Corollary 2.4 and Propositions 2.6, 2.7 and 2.8(a) in \cite{HO},
if a type I or type II edgepath system
is compressible or indeterminate,
then its cycle of final $r$-values 
must include both positive and negative.
Since 
all the signed final $r$-values
of the monotonic edgepath system have the common sign,
the edgepath system is incompressible.
\\
(2): It follows from Proposition 2.7(3)(a) in \cite{HO}. \\
(3)(b): It follows from Proposition 2.9 in \cite{HO}, where the condition (*) is stated. \\
(4): It follows from Proposition 2.5 in \cite{HO}.
\end{proof}

\section{Calculation and Comparison of the twists} 
\label{Sec:Twist}

This section is devoted to proving Theorem \ref{Thm:Diam:LowerBound:ByCrossing}.
The key is the estimation of the maximal and minimal ``twists" defined as follows. 

The \textit{twist} $\Twist(\Edgepath)$ of 
an edgepath $\Edgepath$ is defined as 
the sum of $-2\,\sigma( \Edge_i )\,| \Edge_i |$ for non-$\infty$-edges $\{ e_i \}$ 
included in $\Edgepath$. 
The \textit{twist} $\Twist(\EdgepathSystem)$ of 
an edgepath system $\EdgepathSystem$ is then defined as 
the sum of $\Twist( \Edgepath_i )$ for edgepaths 
$\Edgepath_i$
included in $\EdgepathSystem=\parestart
\Edgepath_1, \Edgepath_2, \ldots, \Edgepath_\NumTangles \pareend$ 
. 
Then, as described in \cite{HO}, 
the boundary slope of an essential surface 
associated to an edgepath system $\EdgepathSystem$ is 
calculated by $\Twist( \EdgepathSystem )-\Twist(\EdgepathSystem_{\mathrm{Seifert}})$, 
where $\EdgepathSystem_{\mathrm{Seifert}}$ denotes 
the edgepath system $\EdgepathSystem$ to which 
a Seifert surface for $K$ is associated. 
Note that $\infty$-edges and constant edgepaths contribute nothing to the twist.

In the following, we call the maximum among 
the twists of all the incompressible or indeterminate candidate edgepath systems 
for a Montesinos knot $K$ the \textit{maximal twist} for $K$. 
The \textit{minimal twist} for $K$ is defined in the same way.

%
%
\subsection{Estimation of maximal and minimal twists} 
Precisely we establish the following estimation about twists.

\begin{proposition}
\label{Prop:Twist:LowerBound}
Let $\Twist_{\max}$ and $\Twist_{\min}$
denote the maximal and minimal twists 
for a Montesinos knot $K$, respectively. 
Let $\BasicEdgepathSystem_\inc$ and $\BasicEdgepathSystem_\dec$
be the monotonically increasing and decreasing basic edgepaths for $K$ 
with the twists $\Twist_\inc$ and $\Twist_\dec$, respectively. 
Then $\Twist_{\max}$ satisfies
\[
 \left\{
  \begin{array}{ll}
   \Twist_{\max}=\Twist_{\dec}+2\,\BasicEdgepathSystem_\dec(0)\ge \Twist_{\dec} 
    & \textrm{if $\BasicEdgepathSystem_\dec(0)\ge 0$,} \\
   \Twist_{\max}\ge\Twist_{\dec}-6 
    & \textrm{if $\BasicEdgepathSystem_\dec(0)= -1$,} \\
   \Twist_{\max}=\Twist_{\dec}
    & \textrm{if $\BasicEdgepathSystem_\dec(0)\le -2$.} \\
  \end{array}
 \right.
\]
Also $\Twist_{\min}$ satisfies
\[
 \left\{
  \begin{array}{ll}
   \Twist_{\min}=\Twist_{\inc}+2\,\BasicEdgepathSystem_\inc(0)\le \Twist_{\inc}
    & \textrm{if $\BasicEdgepathSystem_\inc(0)\le 0$,} \\
   \Twist_{\min}\le \Twist_{\inc} + 6 
    & \textrm{if $\BasicEdgepathSystem_\inc(0)= +1$,} \\
   \Twist_{\min}=\Twist_{\inc}
    & \textrm{if $\BasicEdgepathSystem_\inc(0)\ge +2$.} \\
  \end{array}
 \right.
\]
\end{proposition}

Once the above proposition is established,
Theorem \ref{Thm:Diam:LowerBound:ByCrossing} is proved as follows.

\begin{proof}[Proof of Theorem \ref{Thm:Diam:LowerBound:ByCrossing}]
As remarked in Subsection \ref{Subsec:HO:Montesinos knot}, 
Theorem \ref{Thm:Diam:LowerBound:ByCrossing} holds 
if the number of tangles $N$ in a Montesinos knot $K$ is at most two. 
Thus we assume that $N \ge 3$. 
In the following we use the same notations as in Proposition \ref{Prop:Twist:LowerBound}. 
If $\BasicEdgepathSystem_{\dec}(0)\ge 0$ or $\BasicEdgepathSystem_{\inc}(0)\le 0$ holds,
then $K$ is alternating, and we already have $\Diam(K)=2\,\Crossing(K)$ in \cite{IM2}.
In other cases, as in \cite{IM2}, 
$2\,\Crossing(K)=\Twist_\dec-\Twist_\inc$ holds.
If $\BasicEdgepathSystem_\dec(0)=-1$,
we have $\Twist_{\max}\ge \Twist_\dec-6$ and $\Twist_{\min}= \Twist_\inc$.
If $\BasicEdgepathSystem_\inc(0)=+1$,
we have $\Twist_{\max}= \Twist_\dec$ and $\Twist_{\min}\le \Twist_\inc+6$.
Otherwise,
we have $\Twist_{\max}\ge \Twist_\dec$ and $\Twist_{\min}\le \Twist_\inc$.
Hence, in all cases, 
$\Diam(K)=\Twist_{\max}-\Twist_{\min}\ge \Twist_{\dec}-\Twist_{\inc}-6 =2\,\Crossing(K)-6$.
Note that 
$\BasicEdgepathSystem_\dec(0)=-1$ and $\BasicEdgepathSystem_\inc(0)=+1$
cannot occur at the same time 
since 
the relation $\BasicEdgepathSystem_\inc(0)=\BasicEdgepathSystem_\dec(0)+\NumTangles$ holds 
and we are assuming $\NumTangles\ge 3$. 
\end{proof}

%
%
\subsection{Estimation of twists}
\label{Subsec:CompTwist}
In order to prove Proposition \ref{Prop:Twist:LowerBound}, 
we have to estimate 
the maximal and minimal twists for a Montesinos knot $K$. 
In this subsection,
we prepare two lemmas 
giving formulae for comparing twists of edgepath systems.

\subsubsection{Twists of type I edgepath systems}

Here we give a lemma used to compare the twists of type I edgepath systems.
In the proof, we introduce an integration formula to compute twists of edgepath systems, 
which is of interest independently. 

In the following, we regard an edgepath system 
as a continuous and piecewise-affine function 
$[0,1]\rightarrow\mathbb{R}$ 
as explained in \ref{SubsubSec:Basic:ExtendedBasic:EdgepathSystem}. 
Precisely, for 
a basic edgepath system or a type I edgepath system $\EdgepathSystem$, 
we define 
a continuous and piecewise-affine function $\func:[0,1]\rightarrow\mathbb{R}$ as follows.
If $\EdgepathSystem$ is a basic edgepath system $\BasicEdgepathSystem$,
we define $\func$ as $\widetilde{\BasicEdgepathSystem}$ regarded as a function.
If $\EdgepathSystem$ is of type I,
assuming that $\EdgepathSystem$ is constructed by cutting 
a basic edgepath system $\BasicEdgepathSystem$ at $u=u_0$, 
we define $\func$ so that
$\func(u)=0=\widetilde{\BasicEdgepathSystem}(u_0)$ for $0\le u\le u_0$ and
$\func(u)=\widetilde{\BasicEdgepathSystem}(u)$ for $u_0\le u\le 1$.

\begin{lemma}
\label{Lem:Twist:TypeI}
Let $\BasicEdgepathSystem_\dec$ 
be the monotonically decreasing basic edgepath system for a Montesinos knot $K$. 
\begin{enumerate}
 \item \label{Lem:Twist:TypeI:General}
  For any type I edgepath system $\EdgepathSystem$,
  its twist $\Twist(\EdgepathSystem)$ satisfies
  $\Twist(\EdgepathSystem)
    \le \Twist(\BasicEdgepathSystem_\dec) + 2\BasicEdgepathSystem_{\dec}(0)
    \le \Twist(\BasicEdgepathSystem_\dec) -2$.
 \item \label{Lem:Twist:TypeI:Case2}
  Assume that $\BasicEdgepathSystem_\dec(0)=-1$.
  \begin{itemize}

   \item[\textup{(a)}] 
          Assume further that there exists a solution $u=u_0$ 
          for the equation $\widetilde{\BasicEdgepathSystem_\dec}(u)=0$. 
          Let $\Twist$ be the twist of 
          the type I edgepath system $\EdgepathSystem_{\typeI,\dec}$ 
          obtained by cutting $\BasicEdgepathSystem_\dec$ at $u=u_0$. 
            		\begin{itemize}
			\item[(a1)] \label{Lem:Twist:TypeI:Case2:1/2}
		If $0< u_0\le 1/2$ holds, then $\tau$ satisfies 
          $\Twist(\BasicEdgepathSystem_\dec)-4 \le \Twist < \Twist(\BasicEdgepathSystem_\dec)-2$.
          This twist is maximal among all type I edgepath systems.

		   \item[\textup{(a2)}] \label{Lem:Twist:TypeI:Case2:2/3}
		If $1/2 < u_0\le 2/3$ holds, then $\tau$ satisfies 
          $\Twist(\BasicEdgepathSystem_\dec)-6 \le \Twist < \Twist(\BasicEdgepathSystem_\dec)-4$.
          This twist is maximal among all type I edgepath systems.
			\end{itemize}
			
   \item[\textup(b)] \label{Lem:Twist:TypeI:Case2:Other}
          Assume further that $\BasicEdgepathSystem_\dec(1/2)<0$ holds.
	  Then, for any type I edgepath system $\EdgepathSystem$, 
	  its twist $\Twist (\EdgepathSystem)$ satisfies $\Twist \le \Twist(\BasicEdgepathSystem_\dec)-4$.

  \end{itemize}

 \end{enumerate}
\end{lemma}

\begin{proof}
We first prepare the following claim.

\begin{claim*}
Let $\EdgepathSystem_a$ and $\EdgepathSystem_b$
be a basic or type I edgepath systems.
Assume that $\func_a$ and $\func_b$
are the corresponding functions $[0,1]\rightarrow\mathbb{R}$
for the edgepath systems.
If $\func_a(u)\ge \func_b(u)$ holds for $0\le u\le 1$,
then, their twists satisfy 
$\Twist(\EdgepathSystem_a)\le \Twist(\EdgepathSystem_b)$.
\end{claim*}

\begin{proof}

The following subclaim gives an integration formula to compute twists.

\begin{subclaim}
Let $\EdgepathSystem$ be 
a basic edgepath system or an edgepath system of type I. 
Then its twist is calculated by the integration of the form
\begin{eqnarray}
\Twist(\EdgepathSystem)
&=&
\int_1^{0} 
  -\frac{2}{(u-1)^2} 
  \frac{d \func(u)}{d u}
 du. 
 \label{Eq:TwistByIntegration}
\end{eqnarray}
\end{subclaim}

\begin{proof}
We fix a basic edgepath system $\BasicEdgepathSystem$ and 
the extended basic edgepath system 
$\widetilde{\BasicEdgepathSystem} = \parestart $%
$\widetilde{\BasicEdgepath_{1}}$,
$\widetilde{\BasicEdgepath_{2}}$, \ldots,
$\widetilde{\BasicEdgepath_{\NumTangles}}$%
$\pareend$ for $\BasicEdgepathSystem$. 
Let $\EdgepathSystem_{u_0}=\parestart $%
$\Edgepath_{1,u_0}$,
$\Edgepath_{2,u_0}$, \ldots,
$\Edgepath_{\NumTangles,u_0}$%
$\pareend$ denote
a type I edgepath system obtained from $\BasicEdgepathSystem$ 
by cutting at $u=u_0$ in our manner.
Though $\EdgepathSystem_{u_0}$ may not satisfy the gluing consistency
(\ref{Eq:EquationGluingConsistency}),
we can calculate the twist
$\Twist(\Edgepath_{i,u_0})$ and $\Twist(\EdgepathSystem_{u_0})$ formally.
We regard $\Twist(\EdgepathSystem_{u})$ as a function of $u$ from $[0,1]$ to $\mathbb{R}$.

Assume first that 
$\Edgepath_i$ 
ends
at a point on a decreasing edge $\angleb{p/q}$\,--\,$\angleb{r/s}$,
where $ps-qr=-1$.
Let $\alpha$ denote the twist of complete edges of $\Edgepath_i$
included completely in $u\ge u_0$.
The derivative of the twist of the edgepath $\Edgepath_{u_0,i}$ 
at $u_0$ is calculated as
\[
\frac{d \Twist( \Edgepath_{i,u_0} )}{d u_0}
=\frac{d }{d u_0}
 \left(
  \alpha+2\,\frac{1+s(u_0-1)}{(s-q)(u_0-1)}
 \right)
=-\frac{2}{ (s-q)(u_0-1)^2}
\]
with the formula (\ref{Eq:Formula:LengthOfPartialEdge}) of length of a partial edge 
in Subsubsection \ref{SubsubSec:LengthOfEdgepath}.
Note that $\alpha$ is constant at $(q-1)/q\le u_0 \le (s-1)/s$.
Besides,
\begin{eqnarray*}
\frac{d \widetilde{\BasicEdgepath_{i}} (u_0)}{d u_0} 
&=&
 \frac{d}{du_0}
 \left( 
 \frac{r}{s}+\left(\frac{p}{q}-\frac{r}{s}\right)
  \frac{ \frac{s-1}{s}-u_0 }{ \frac{s-1}{s} - \frac{q-1}{q} }
 \right)
= \frac{1}{s-q}
.
\end{eqnarray*}
Hence,
\begin{eqnarray}
\label{Eq:dTwistdu}
\frac{d \Twist( \Edgepath_{i,u_0} )}{d u_0}
&=& -\frac{2}{(u_0-1)^2} \frac{d \widetilde{\BasicEdgepath_{i}} (u_0)}{d u_0}
.
\end{eqnarray}
This identity also holds if $\Edgepath_i$ 
ends
at a point on an increasing edge or a horizontal edge.
Since $\Twist(\EdgepathSystem_{u_0})$ is $0$ for $u_0$ close to $1$,
by summing up (\ref{Eq:dTwistdu}) and performing integration,
we have
\begin{eqnarray*}
\Twist(\EdgepathSystem_{u_0})
&=&
\Twist(\EdgepathSystem_{u_0})
-\Twist(\EdgepathSystem_{1})
=
\sum_{i=1}^{\NumTangles}
 (\Twist(\Edgepath_{i,u_0})
 -\Twist(\Edgepath_{i,1}))
=
\sum_{i=1}^{\NumTangles}
\int_1^{u_0} 
  \frac{d \Twist(\Edgepath_{i,u})}{d u}
 du 
\\
&=&
\sum_{i=1}^{\NumTangles}
\int_1^{u_0} 
 \left(
  -\frac{2}{(u-1)^2} \frac{d \widetilde{\BasicEdgepath_{i}} (u)}{d u}
 \right)
 du 
=
\int_1^{u_0} 
  -\frac{2}{(u-1)^2} 
  \frac{d \widetilde{\BasicEdgepathSystem}(u) }{d u}
 du 
.
\\
&=&
\int_1^{0} 
  -\frac{2}{(u-1)^2} 
  \frac{d \func(u)}{d u}
 du. 
\end{eqnarray*}
\end{proof}

By this subclaim, we have; 

\begin{eqnarray*}
\Twist(\EdgepathSystem_a)-\Twist(\EdgepathSystem_b)
&=&
\int_1^{0} 
  -\frac{2}{(u-1)^2} 
  \frac{d \func_a(u) }{d u}
 du
 - 
\int_1^{0} 
  -\frac{2}{(u-1)^2} 
  \frac{d \func_b(u) }{d u}
 du
\\
&=&
\int_1^{0} 
  -\frac{2}{(u-1)^2} 
  \frac{d (\func_a(u)-\func_b(u)) }{d u}
 du
\\
&=&
\left[
  -\frac{2}{(u-1)^2} 
  ( \func_a(u)-\func_b(u) )
\right]_{1}^{0}
\\
&&-
\int_1^{0} 
  \frac{d}{du} \left(-\frac{2}{(u-1)^2} \right)
  \cdot
  ( \func_a(u)-\func_b(u) )
 du
\\
&\le& 0
.
\end{eqnarray*}
Note that $\func_a(0)=\func_b(0)=0$ holds by the gluing consistency.
Besides, $\EdgepathSystem_a$ and $\EdgepathSystem_b$ comes from the same Montesinos knot $K=M(T_1,T_2,\ldots,T_\NumTangles)$.
Edges of the edgepaths $\Edgepath_{a,i}$ and $\Edgepath_{b,i}$ near $u=1$ are both the horizontal edge corresponding to $T_i$ .
Hence, there exists $u_1<1$ such that
$\func_a(u)=\func_b(u)$ holds for $u_1<u\le 1$.
Eventually, the square bracket has the value zero.
%
\end{proof}


In the following, 
let $\func_\dec$ denote
the function 
for the monotonically decreasing basic edgepath system
$\BasicEdgepathSystem_\dec$. \\

(1): 
Let $\func_1:[0,1]\rightarrow \mathbb{R}$ be a type I edgepath system $\EdgepathSystem$ as a function.
From $\BasicEdgepathSystem_\dec$,
by replacing the last segment $(0,z)$\,--\,$(u_1,z+v_1)$
with a polygonal line $(0,0)$\,--\,$(\varepsilon,z)$\,--\,$(u_1,z+v_1)$,
we have a function $\func_2:[0,1]\rightarrow \mathbb{R}$,
which is pretty close to $\BasicEdgepathSystem_\dec$ in $[\varepsilon,1]$
and is pretty close to vertical edges connecting
 $\angleb{\BasicEdgepathSystem_\dec(0)}$ and $\angleb{0}$.
Since $\func_1(u) \ge \func_2(u)$ holds for $0\le u\le 1$,
by the method in the above claim,
we have
\begin{eqnarray*}
 \Twist(\EdgepathSystem)
 &=&
  \int_{1}^{0} -\frac{2}{(u-1)^2} \frac{d\varphi_1 (u)}{du} du
 \le  
  \int_{1}^{0} -\frac{2}{(u-1)^2} \frac{d\varphi_2 (u)}{du} du
\end{eqnarray*}
for each sufficiently small $\varepsilon>0$.
Now,
\begin{eqnarray*}
 \lim_{\varepsilon\rightarrow 0}
  \int_{1}^{0} -\frac{2}{(u-1)^2} \frac{d\varphi_2 (u)}{du} du
 &=&
 \lim_{\varepsilon\rightarrow 0}
  \left(
   \int_\varepsilon^{0} -2 \frac{d\varphi_2 (u)}{du} du
   +\int_{1}^{\varepsilon} -\frac{2}{(u-1)^2} \frac{d\varphi_2 (u)}{du} du
  \right)
 \\
 &=&
  2 \BasicEdgepathSystem_\dec(0)
  + \Twist(\BasicEdgepathSystem_\dec)
  .
\end{eqnarray*}
Hence,
\begin{eqnarray*}
 \Twist(\EdgepathSystem)
 &\le&
  2 \BasicEdgepathSystem_\dec(0)
  + \Twist(\BasicEdgepathSystem_\dec)
 \le
  \Twist(\BasicEdgepathSystem_\dec) -2
.
\end{eqnarray*}
Note that $\BasicEdgepathSystem_\dec(0)\le -1$
 is necessary to give a type I edgepath system.


 (\ref{Lem:Twist:TypeI:Case2:1/2})(a1):   
  Let $\func_{\typeI,\dec}$ denote the function corresponding to 
  the monotonically decreasing type I edgepath system $\EdgepathSystem_{\typeI,\dec}$. 
  Then $\func_{\typeI,\dec}$ and $\func_\dec$ coincide for $u\ge u_0$,
  while
  $\func_{\typeI,\dec}(u)=0$ and $\func_\dec(u)=u/u_0-1$ hold for $0\le u\le u_0$.
  By the above subclaim,
  \begin{eqnarray*}
    \Twist(\EdgepathSystem_{\typeI,\dec})-\Twist(\BasicEdgepathSystem_\dec)
    &=&
    -\int_{u_0}^{0}
     -\frac{2}{(u-1)^2} 
     \frac{d \func_\dec(u) }{d u} 
    du
    \\
    &=&
    -\int_{u_0}^{0}
     -\frac{2}{(u-1)^2} 
     \frac{1}{u_0} 
    du
    = -\frac{2}{1-u_0} \in [-4,-2) 
    .
  \end{eqnarray*}

  Suppose that there exists another type I edgepath system $\EdgepathSystem$ 
  whose twist is greater than $\Twist(\EdgepathSystem_{\typeI,\dec})$ and 
  which is obtained by cutting a basic edgepath system $\BasicEdgepathSystem$ at $u = u_1$.
  Then we note that $\BasicEdgepathSystem$ 
  must include a part below the horizontal line $v=0$ by the above claim.
  On the other hand, $\BasicEdgepathSystem$ satisfies 
  $\BasicEdgepathSystem(u_1)=0$,
  $\BasicEdgepathSystem(u_2)<0$ and
  $\BasicEdgepathSystem(1/2)\ge 0$ for some $u_1$ and $u_2$ 
  satisfying $u_1<u_2<u_0\le 1/2$.
  Though, this is impossible since $\BasicEdgepathSystem$ is
  affine in the interval $0\le u\le 1/2$.
  %


 (\ref{Lem:Twist:TypeI:Case2:2/3})(a2): 
  Again,  let $\func_{\typeI,\dec}$ denote the function corresponding to 
  the monotonically decreasing type I edgepath system $\EdgepathSystem_{\typeI,\dec}$. 
  Let $\func_1$ be a function defined by 
  $\func_1 (u) = u/u_0-1$ for $0\le u\le u_0$ and 
  $\func_1 (u) = \func_\dec(u)$ for $u_0\le u\le 1$. 
  The graph of $\func_1$ includes a segment connecting $(0,-1)$ and $(u_0,0)$. 
  Thus $\func_{\dec}(u) \ge \func_1(u)$ holds for $0\le u\le 1$.
  Note that $\BasicEdgepathSystem_\dec$ is a concave function.
  Then, by the above claim and subclaim 
  we have
  \begin{eqnarray*}
    \Twist(\EdgepathSystem_{\typeI,\dec})-\Twist(\BasicEdgepathSystem_\dec)
    &\ge& \Twist(\EdgepathSystem_{\typeI,\dec})-\Twist(\func_1)
    \\
    &=& -\int_{u_0}^{0}
     -\frac{2}{(u-1)^2} 
     \frac{d \func_1(u) }{d u} 
    du
    =
    -\int_{u_0}^{0}
     -\frac{2}{(u-1)^2} 
     \frac{1}{u_0} 
    du
    \\
    &=& -\frac{2}{1-u_0}\ge -6
    .
  \end{eqnarray*}
  Also let $\func_2$ be a function defined by 
  $\func_2(u)=2u-1$ for $0\le u\le 1/2$, 
  $\func_2(u)=0$ for $1/2 \le u \le u_0$, and 
  $\func_2(u)=\func_\dec(u)$ for $u_0\le u\le 1$.
  The graph of $\func_2$ includes a piecewise affine path 
  connecting $(0,-1)$, $(1/2,0)$ and $(u_0,0)$.
  Thus $\func_{\typeI,\dec}(u) \ge \func_2(u)$ holds for $0\le u\le 1$.
  Again, by the above claim and subclaim, we have
  \begin{eqnarray*} 
    \Twist(\EdgepathSystem_{\typeI,\dec})-\Twist(\BasicEdgepathSystem_\dec)
    &<& \Twist(\EdgepathSystem_{\typeI,\dec})-\Twist(\func_2)
    \\
    &=& -\int_{\frac{1}{2}}^{0}
     -\frac{2}{(u-1)^2} 
     \frac{d \func_2(u) }{d u} 
    du
    =
    -\int_{\frac{1}{2}}^{0}
     -\frac{2}{(u-1)^2} 
     \cdot 2
    du
    \\
    &=& -4
    .
  \end{eqnarray*}
  Here, $\Twist(\func)$ denotes the integration value of
 (\ref{Eq:TwistByIntegration}) for the function $\func$.

  Suppose that $\EdgepathSystem$ is another type I edgepath system 
  whose twist $\Twist$ is greater than $\Twist(\EdgepathSystem_{\typeI,\dec})$.
  A basic edgepath system $\BasicEdgepathSystem$ including $\EdgepathSystem$
  must include a part below the line $v=0$.
  Since $\BasicEdgepathSystem$ is affine in the intervals
  $0<u<1/2$ and $1/2<u<2/3$,
  $\BasicEdgepathSystem$ satisfies 
  $\BasicEdgepathSystem(0)>0$, $\BasicEdgepathSystem(1/2)<0$
  and $\BasicEdgepathSystem(2/3)\ge0$.
  Though, this is impossible as follows.
  $\BasicEdgepathSystem(1/2)<0$ means that
  only one of the edgepaths satisfies
  $z+1/2 \le \BasicEdgepath_i(1/2)<z+1$ for some integer $z$,
  and thus,
  at most one of the basic edgepaths can end with an increasing edge.
  Thus, $\BasicEdgepathSystem(0)=0$ or $-1$.
  %


 (\ref{Lem:Twist:TypeI:Case2:Other})(b)
  Assume that $\EdgepathSystem$ is
  a type I edgepath system with the function $\func$.
  Let $\func_1$ be a function defined by 
  $\func_1(u)=2\cdot\BasicEdgepathSystem_\dec(1/2)\cdot u$ for $0\le u \le 1/2$ and 
  $\func_1(u)= \BasicEdgepathSystem_\dec(u)$ for $1/2\le u\le 1$.
  The graph of $\func_1$ includes a segment
  connecting $(0,0)$ and $(1/2,\BasicEdgepathSystem_\dec(1/2))$.
  Then, $\func(u)\ge \func_1(u)\ge \func_\dec(u)$ holds for $0\le u\le 1$.
  Note that
  $\func_\dec(u)=2\cdot(\BasicEdgepathSystem_\dec(1/2)+1)u-1$ holds for $0\le u \le 1/2$.
  By the above claim and subclaim, we have;
  \begin{eqnarray*}
    \Twist(\EdgepathSystem)-\Twist(\BasicEdgepathSystem_\dec)
    &\le&
    \Twist(\func_1)-\Twist(\BasicEdgepathSystem_\dec)
    \\
    &=& \int_{1}^{0}
     -\frac{2}{(u-1)^2} 
     \frac{d (\func_1(u)-\func_\dec(u)) }{d u} 
    du
    =
    \int_{\frac{1}{2}}^{0}
     -\frac{2}{(u-1)^2} 
     \cdot (-2)
    du
    \\
    &=& -4
    .
  \end{eqnarray*}

\end{proof}

\subsubsection{Twists of type II or III edgepath systems}

The second lemma is also for calculating twist
and is useful for type II and type III edgepath systems.
We first introduce a classification of edgepath systems.

\begin{definition}\ \\ \vspace{-4mm}
 \begin{enumerate}
  \item For a basic edgepath $\BasicEdgepath$, we call it:
    {\rm (a)} 
           a class A basic edgepath if it is a monotonically decreasing basic edgepath $\BasicEdgepath_{\dec}$, 
    {\rm (b)} 
           a class B basic edgepath if $\BasicEdgepath$ and $\BasicEdgepath_{\dec}$ bound 
           single triangle and single vertical edge in the strip $\Strip$, 
    {\rm (c)} 
           a class C basic edgepath if $\BasicEdgepath$ and $\BasicEdgepath_{\dec}$ bound 
           two triangles and no vertical edges in the strip $\Strip$. 
    A class B basic edgepath is obtained
    from the monotonically decreasing basic edgepath
    by replacing an edge $\angleb{z+0}$\,--\,$\angleb{z+1/2}$ 
    by an edge $\angleb{z+1}$\,--\,$\angleb{z+1/2}$ for some integer $z$.
  \item For a basic edgepath system $\BasicEdgepathSystem$, we call it:
    {\rm (a)} 
          a class A basic edgepath system
          if it is a monotonically decreasing basic edgepath system,
    {\rm (b)} 
          a class B/C basic edgepath system 
          if exactly one basic edgepath in the basic edgepath system is class B/C
          and all the other edgepaths are monotonically decreasing (class A),
  \item For a type II edgepath system,
        it is called class A/B/C
        if it is obtained from a class A/B/C basic edgepath system
        by extending by vertical edges.
  \item For a type III edgepath system,
        it is called class A/B/C
        if it is obtained from a class A/B/C basic edgepath system
        by extending by $\infty$-edges for all edgepaths.
 \end{enumerate}
\end{definition}

%
%

  \begin{figure}[htb]
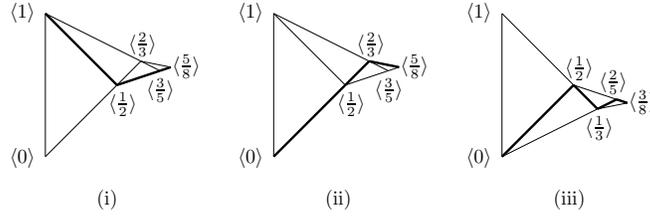


   \begin{minipage}{83pt}
   \begin{center}
    \scalebox{0.75}{
    \begin{picture}(80,100)
     \put(0,20){\scalebox{0.4}{\includegraphics{classb.eps}}}
     \put(30,0){(i)}
     \put(-14,21){$\angleb{0}$}
     \put(-14,94){$\angleb{1}$}
     \put(36,48){$\angleb{\frac{1}{2}}$}
     \put(46,78){$\angleb{\frac{2}{3}}$}
     \put(55,55){$\angleb{\frac{3}{5}}$}
     \put(68,67){$\angleb{\frac{5}{8}}$}
    \end{picture}
    }
   \end{center}
   \end{minipage}
    \begin{minipage}{83pt}
    \begin{center}    
    \scalebox{0.75}{
     \begin{picture}(80,100)
      \put(0,20){\scalebox{0.4}{\includegraphics{classc.eps}}}
      \put(30,0){(ii)}
      \put(-14,21){$\angleb{0}$}
      \put(-14,94){$\angleb{1}$}
      \put(36,48){$\angleb{\frac{1}{2}}$}
      \put(46,78){$\angleb{\frac{2}{3}}$}
      \put(55,55){$\angleb{\frac{3}{5}}$}
      \put(68,67){$\angleb{\frac{5}{8}}$}
     \end{picture}
    }
    \end{center}
    \end{minipage}
    \begin{minipage}{83pt}
    \begin{center}    
    \scalebox{0.75}{
     \begin{picture}(80,100)
      \put(0,20){\scalebox{0.4}{\includegraphics{classc2.eps}}}
      \put(30,0){(iii)}
      \put(-14,21){$\angleb{0}$}
      \put(-14,94){$\angleb{1}$}
      \put(36,66){$\angleb{\frac{1}{2}}$}
      \put(46,36){$\angleb{\frac{1}{3}}$}
      \put(53,59){$\angleb{\frac{2}{5}}$}
      \put(68,49){$\angleb{\frac{3}{8}}$}
     \end{picture}
    }
    \end{center}
    \end{minipage}

    \caption{
    (i) A class B basic edgepath for $\angleb{5/8}$,
    (ii) a class C basic edgepath for $\angleb{5/8}$ and
    (iii) a not-class-C basic edgepath for $\angleb{3/8}$.
    Note that the edgepath in (iii) 
    is not class C since it is not minimal.
    }
    \label{Fig:Edgepaths}

  \end{figure}

Then we obtain the following lemma.

%
\begin{lemma} 
\label{Lem:Twist:TypeIIorIII}
Let $\BasicEdgepathSystem_\dec$ be the monotonically decreasing edgepath system 
with twist $\Twist(\BasicEdgepathSystem_\dec )$ for a Montesinos knot $K$. 
 \begin{enumerate}
  \item 
	Let $\EdgepathSystem_{\typeII}$ be a type II edgepath system 
	with twist $\Twist (\EdgepathSystem_{\typeII}) $. 
	Then, 
        \label{Lem:Twist:TypeIIorIII:TypeII}
    {\rm (a)} 
          $\Twist (\EdgepathSystem_{\typeII}) = \Twist (\BasicEdgepathSystem_\dec) + 2\,\BasicEdgepathSystem_\dec(0)$ 
          if and only if 
          $\EdgepathSystem_{\typeII}$ is of class A,
    {\rm (b)} 
          $\Twist (\EdgepathSystem_{\typeII}) = \Twist (\BasicEdgepathSystem_\dec) + 2\,\BasicEdgepathSystem_\dec(0) -2 $ 
          if and only if 
          $\EdgepathSystem_{\typeII}$ is of class B,
    {\rm (c)} 
          otherwise, 
          $\Twist ( \EdgepathSystem_{\typeII} ) \le \Twist (\BasicEdgepathSystem_\dec) +2\,\BasicEdgepathSystem_\dec(0)-4$.
  \item 
	Let $\EdgepathSystem_{\typeIII}$ be a type III edgepath system 
	with twist $\Twist (\EdgepathSystem_{\typeIII}) $. 
	Then, 
        \label{Lem:Twist:TypeIIorIII:TypeIII}
    {\rm (a)} 
          $\Twist (\EdgepathSystem_{\typeIII}) = \Twist (\BasicEdgepathSystem_\dec)$
          if and only if 
          $\EdgepathSystem_{\typeIII}$ is monotonically decreasing, 
    {\rm (b)} 
          $\Twist (\EdgepathSystem_{\typeIII}) = \Twist (\BasicEdgepathSystem_\dec)-4$
          if and only if 
          $\EdgepathSystem_{\typeIII}$ is of class B or C,
    {\rm (c)} 
          otherwise, 
          $\Twist (\EdgepathSystem_{\typeIII}) \le \Twist (\BasicEdgepathSystem_\dec)-6$.
 \end{enumerate}
\end{lemma}
\begin{proof}
Let $\BasicEdgepathSystem_\dec=\parestart \BasicEdgepath_{\dec,1}, \BasicEdgepath_{\dec,2}, \ldots, \BasicEdgepath_{\dec,\NumTangles}\pareend$ be the monotonically decreasing basic edgepath system.

We first prepare two combinatorial values defined for edgepath systems. 
We consider 
a basic edgepath system $\BasicEdgepathSystem$$=\parestart $$\BasicEdgepath_1$, $\BasicEdgepath_2$, $\ldots$, $\BasicEdgepath_\NumTangles$$\pareend$ for $K$.
Let $L_i$ and $V_i$ denote
the number of triangles and the number of vertical edges in the strip $\Strip$
bounded by $\BasicEdgepath_{i}$ and $\BasicEdgepath_{\dec,i}$
 respectively.
Set $L=\sum_{i=1}^\NumTangles L_i$ and $V=\sum_{i=1}^\NumTangles V_i$.

Then we have the following claim. 

\begin{claim*} 
Let $L$ and $V$ denote the corresponding values for $\BasicEdgepathSystem$.
Then the twist $\Twist (\BasicEdgepathSystem)$ is calculated by
\[
\Twist(\BasicEdgepathSystem)= \Twist(\BasicEdgepathSystem_\dec ) -2\,(L+V)
.
\]
\end{claim*}

\begin{proof}
Let $\{t_1, t_2, \ldots, t_L\}$ be the triangles bounded by 
$\BasicEdgepath_{i}$'s and $\BasicEdgepath_{\dec,i}$'s. 
Let $\vec{l}_i$ be the directed loop bounding the triangle $t_i$ 
with anti-clockwise direction.
Let $\{\vec{v}_1, \vec{v}_2, \ldots, \vec{v}_V\}$ be the bounded downward vertical edges.
Let $\overrightarrow{\BasicEdgepathSystem_\dec}$ and $\overrightarrow{\BasicEdgepathSystem}$ denote the sets of
directed edgepaths of $\BasicEdgepathSystem_\dec$ and $\BasicEdgepathSystem$
with the right-to-left direction.
For a leftward edge $\vec{\Edge}$ and 
$-\vec{\Edge}$ with the opposite direction,
we here define
the effect of $-\vec{\Edge}$ on twist
to be the negative of the effect of $\vec{\Edge}$.
The effects of $\vec{\Edge}$ and $-\vec{\Edge}$ on twist cancel out.
Then,
directed arcs and loops satisfy the following identity.
\[
\sum_{\vec{e}\in \overrightarrow{\BasicEdgepathSystem_\dec} \cup (\cup_{j=1}^{L} \vec{l}_j)} \Twist(\vec{e})
=
\sum_{\vec{e}\in \overrightarrow{\BasicEdgepathSystem} \cup (\cup_{j=1}^{V} \vec{v}_j)} \Twist(\vec{e})
.
\]
Now, three vertices of
a triangle in the strip $\Strip$
can be expressed as
$P_1=\angleb{p/q}=((q-1)/q,p/q)$, 
$P_2=\angleb{r/s}=((s-1)/s,r/s)$
and $P_3=\angleb{(p+r)/(q+s)}=((q+s-1)/(q+s),(p+r)/(q+s))$.
Note that 
the $u$-coordinate of $P_3$ is greater than those of $P_1$ and $P_2$,
and
the $v$-coordinate of $P_3$ lies between those of $P_1$ and $P_2$.
Thus, the shapes of triangles are somehow restricted.
As the result, 
contribution of $\vec{l}_j$ is $-2$.
Besides, contribution of a downward vertical edge $\vec{v}_j$ is $+2$.
Thus,
we have $\Twist(\BasicEdgepathSystem_\dec) - 2 L = \Twist(\BasicEdgepathSystem) + 2 V$.
The twist of $\BasicEdgepathSystem$ is calculated by
$\Twist(\BasicEdgepathSystem)=\Twist(\BasicEdgepathSystem_\dec)-2\,(L+V)$.
\end{proof}

(\ref{Lem:Twist:TypeIIorIII:TypeII})
A type II edgepath system $\EdgepathSystem_{\typeII}$ is obtained from
a basic edgepath system $\BasicEdgepathSystem$ by extending by vertical edges.
Let $L$ and $V$ denote the corresponding values for $\BasicEdgepathSystem$. 
Now, the effect on the twist of $\EdgepathSystem_{\typeII}$ by vertical edges
is $2\,\BasicEdgepathSystem_\dec(0)+2V$.
Thus, by the claim above, we have; 
$\Twist(\EdgepathSystem_{\typeII})
=\Twist(\BasicEdgepathSystem) +2\BasicEdgepathSystem_\dec(0)+2V
=\Twist(\BasicEdgepathSystem_\dec) -2\,(L+V)+2\BasicEdgepathSystem_\dec(0)+2V
= \Twist(\BasicEdgepathSystem_\dec) -2L+2\BasicEdgepathSystem_\dec(0)$.
Note that if $L\le 1$, then $V \le 1$ holds, 
and $(L,V)=(1,0)$ is non-minimal since $(L_i,V_i)=(1,0)$ must hold for some $i$. 
Thus, if $L\le 1$, then possible $(L,V)$ is $(0,0)$ or $(1,1)$.
These correspond to class A and class B type II edgepath systems.

(\ref{Lem:Twist:TypeIIorIII:TypeIII})
The twist of a type III edgepath system
is same as its corresponding basic edgepath system $\BasicEdgepathSystem$.
Let $L$ and $V$ denote the corresponding values for $\BasicEdgepathSystem$. 
Then, by the claim above, 
the twist $\Twist (\BasicEdgepathSystem)$ is calculated as 
$\Twist (\BasicEdgepathSystem) =\Twist(\BasicEdgepathSystem_\dec) -2\,(L+V)$.
Here we consider the cases satisfying $L+V\le 2$.
By definition, $V_i$ is $0$ or $1$ for each $i$.
If $V_i=1$, then $L_i\ge 1$,
since some increasing edges must appear in $\BasicEdgepath_i$.
Hence, $L_i\ge V_i$ holds for each $i$, and so does $L\ge V$.
If $(L,V)=(1,0)$, then $(L_i,V_i)=(1,0)$ holds for some $i$,
and this means that $\BasicEdgepath_i$ is not minimal.
Eventually, if $L+V\le 2$,
then $(L,V)$ is $(0,0)$, $(1,1)$ or $(2,0)$.
Moreover, if $(L,V)$ is $(1,1)$ or $(2,0)$,
then $(L_i,V_i)$ is $(1,1)$ or $(2,0)$ for the unique $i$ respectively, 
since $(L_i,V_i)=(1,0)$ is impossible.
Thus, at most one basic edgepath of a basic edgepath system
can be non-monotonically-decreasing.
This implies the assertions (a), (b) and (c) in (\ref{Lem:Twist:TypeIIorIII:TypeIII}).

\end{proof}

%
\subsection{Proof of Proposition \ref{Prop:Twist:LowerBound}} 
\label{Subsec:Proof:Prop:Twist:LowerBound}

This subsection is devoted to 
proving Proposition \ref{Prop:Twist:LowerBound}. 

In the following proof, 
we divide the set of Montesinos knots into some classes
by the nature of their monotonic basic edgepath systems.
For each class, we find an edgepath system 
which is incompressible or indeterminate, 
and then calculate or estimate its twist.
Since the edgepath system is taken
so that its twist gives an appropriate bound of the maximal/minimal twist,
it suffices to prove Proposition \ref{Prop:Twist:LowerBound}.
In fact, our main targets to study are 
monotonic or nearly monotonic edgepath systems. 

\begin{proof}[Proof of Proposition \ref{Prop:Twist:LowerBound}.]

We only prove the assertion about the maximal twist.
It is sufficient 
since the assertion about the minimal twist must hold 
by symmetry.
That is,
with the fact that 
$K=M(\Tangle_1$, $\Tangle_2$, $\ldots,$ $\Tangle_\NumTangles)$
has its mirror image $K^\prime=M(-\Tangle_1$, $-\Tangle_2$, $\ldots,$ $-\Tangle_\NumTangles)$,
the assertion about the minimal twist of the knot $K$
is immediately obtained from 
the assertion about the maximal twist of the knot $K^\prime$.

Let $\BasicEdgepathSystem_\dec$ be the monotonically decreasing edgepath system 
with twist $\Twist(\BasicEdgepathSystem_\dec )$ for a Montesinos knot $K$.

By Lemma \ref{Lem:Incompressibility},
depending on the value of $\BasicEdgepathSystem_\dec(0)$,
we may be able to obtain an incompressible type II or type III edgepath system.
Hence, 
we first divide the Montesinos knots by the value of $\BasicEdgepathSystem_\dec(0)$,
as follows.
\begin{description}
\item[Case 1] 
 $K$ satisfies 
 $\BasicEdgepathSystem_\dec(0) \ge 0$.
\item[Case 2] 
 $K$ satisfies 
 $\BasicEdgepathSystem_\dec(0) = -1$.
\item[Case 3] 
 $K$ satisfies 
 $\BasicEdgepathSystem_\dec(0) \le -2$.
\end{description}
Remember that $\BasicEdgepathSystem_\dec(0)$ is an integer.

\setcounter{internalclaimcounter}{0}


\begin{internalclaim}[Case 1]
For a Montesinos knot with $\BasicEdgepathSystem_\dec(0) \ge 0$,
there exists a monotonically decreasing type II edgepath system $\EdgepathSystem_{\typeII,\dec}$ 
such that an essential surface is associated to $\EdgepathSystem_{\typeII,\dec}$ and 
its twist $\Twist (\EdgepathSystem_{\typeII,\dec})$ satisfies that 
$\Twist (\EdgepathSystem_{\typeII,\dec}) \ge \Twist (\BasicEdgepathSystem_\dec) $. 
\end{internalclaim}
\begin{proof}
From $\BasicEdgepathSystem_\dec$,
by extending one of its edgepaths by some downward vertical edges if necessary,
we have a monotonically decreasing type II edgepath system $\EdgepathSystem_{\typeII, \dec}$ in this case. 
This is incompressible by 
Lemma \ref{Lem:Incompressibility}(\ref{Lem:Incompressibility:MonotonicTypeII})(a).
Its twist $\Twist (\EdgepathSystem_{\typeII,\dec}) $ 
satisfies 
$\Twist (\EdgepathSystem_{\typeII,\dec}) 
=\Twist (\BasicEdgepathSystem_\dec) + 2\,\BasicEdgepathSystem_\dec(0) \ge \Twist (\BasicEdgepathSystem_\dec) $
by Lemma \ref{Lem:Twist:TypeIIorIII}(\ref{Lem:Twist:TypeIIorIII:TypeII})(a).
\end{proof}

\begin{internalclaim}[Case 3]
For a Montesinos knot with $\BasicEdgepathSystem_\dec(0) \le -2$,
the monotonically decreasing type III edgepath system $\EdgepathSystem_{\typeIII,\dec}$ 
such that an essential surface is associated to $\EdgepathSystem_{\typeIII,\dec}$ and 
its twist $\Twist (\EdgepathSystem_{\typeIII,\dec})$ satisfies that 
$\Twist (\EdgepathSystem_{\typeIII,\dec}) = \Twist (\BasicEdgepathSystem_\dec) $. 
\end{internalclaim}
\begin{proof}
Since $\BasicEdgepathSystem_{\dec}(0)\le -2$ in this case, 
$\EdgepathSystem_{\typeIII,\dec}$ is incompressible 
by Lemma \ref{Lem:Incompressibility}(\ref{Lem:Incompressibility:TypeIII}). 
%
Its twist is $\Twist ((\EdgepathSystem_{\typeIII,\dec})) =\Twist(\BasicEdgepathSystem_{\dec})$ by
Lemma \ref{Lem:Twist:TypeIIorIII}(\ref{Lem:Twist:TypeIIorIII:TypeIII})(a).
\end{proof}

As we saw in Subsection \ref{Subsec:HO:Incompressibility}, 
the final $r$-values of a basic edgepath system 
are important in showing the existence of a type II edgepath system with an 
incompressible surface 
by Lemma \ref{Lem:Incompressibility}(\ref{Lem:Incompressibility:MonotonicTypeII})(b).
%
Let $\parestart r_1$, $r_2$, $\ldots,$ $r_\NumTangles$$\pareend$ denote the final $r$-values of the monotonically decreasing basic edgepath system.
Note that all $r_i$'s are negative naturally.
We here divide Case 2 into several cases
according to the numbers of basic edgepaths of $r_i=-1$, $r_i=-2$ or $r_i\le -3$,
as follows.

\begin{description}
 \item[Case 2-1] $K$ satisfies $\sharp\{i\,|\,r_i=-1\}=0$.
 \item[Case 2-2] $K$ satisfies $\sharp\{i\,|\,r_i=-1\}=1$.
  \begin{description}
   \item[Case 2-2-1] Moreover, $K$ satisfies $\sharp\{i\,|\,r_i\le -3\}=0$.
   \item[Case 2-2-2] Moreover, $K$ satisfies $\sharp\{i\,|\,r_i\le -3\}=1$.
    \begin{description}
     \item[Case 2-2-2-1] Moreover, $K$ satisfies $\sharp\{i\,|\,r_i=-2\}=1$. Or equivalently, $\NumTangles=3$.
     \item[Case 2-2-2-2] Moreover, $K$ satisfies $\sharp\{i\,|\,r_i=-2\}\ge 2$. Or equivalently, $\NumTangles\ge 4$.
    \end{description}
   \item[Case 2-2-3] Moreover, $K$ satisfies $\sharp\{i\,|\,r_i\le -3\}\ge 2$.
  \end{description}
 \item[Case 2-3] $K$ satisfies $\sharp\{i\,|\,r_i=-1\}\ge 2$.
  \begin{description}
   \item[Case 2-3-1] Moreover, $\BasicEdgepathSystem_\dec$ of $K$
	      satisfies (*) in 
              Lemma \ref{Lem:Incompressibility}(\ref{Lem:Incompressibility:MonotonicTypeII})(b).
   \item[Case 2-3-2] Moreover, $\BasicEdgepathSystem_\dec$ of $K$ does
	      not satisfy (*) in 
              Lemma \ref{Lem:Incompressibility}(\ref{Lem:Incompressibility:MonotonicTypeII})(b).
  \end{description}
\end{description}
%

\begin{internalclaim}[Cases 2-1, 2-2-3, 2-3-2]
For a Montesinos knot in these cases,
an essential surface is associated to
a type II edgepath system $\EdgepathSystem_{\typeII,\classA}$ obtained by
extending the monotonically decreasing basic edgepath system
by one upward vertical edge.
It has the twist $\Twist=\Twist(\BasicEdgepathSystem_\dec)-2$.
\end{internalclaim}
\begin{proof}
$\BasicEdgepathSystem_\dec$
does not satisfy the condition (*)
in Lemma \ref{Lem:Incompressibility}(\ref{Lem:Incompressibility:MonotonicTypeII})(b)
in all these three cases: 
For there is no $r_i=-1$ in Case 2-1, 
while there is $2$ or more $r_i\le -3$ despite $\sharp\{i\,|\,r_i=-1\}=1$ in Case 2-2-3.
Thus, by Lemma \ref{Lem:Incompressibility}(\ref{Lem:Incompressibility:MonotonicTypeII})(b),
we have a class A type II edgepath system $\EdgepathSystem_{\typeII, \classA}$
to which an incompressible surface is associated.
By Lemma \ref{Lem:Twist:TypeIIorIII}(\ref{Lem:Twist:TypeIIorIII:TypeII})(a),
$\EdgepathSystem_{\typeII,\classA}$ has twist 
$\Twist
=\Twist(\BasicEdgepathSystem_\dec)+2\,\BasicEdgepathSystem_\dec(0)
=\Twist(\BasicEdgepathSystem_\dec)-2$.
\end{proof}

%
\begin{internalclaim}[Cases 2-2-1, 2-2-2-2, 2-3-1]
For a Montesinos knot in these cases,
the monotonically decreasing type I edgepath system 
$\EdgepathSystem_{\typeI,\dec}$ exists,
an essential surface is associated to $\EdgepathSystem_{\typeI,\dec}$, 
and
$\EdgepathSystem_{\typeI,\dec}$ has twist $\Twist$ satisfying 
$\Twist(\BasicEdgepathSystem_\dec)-4\le\Twist<\Twist(\BasicEdgepathSystem_\dec)-2$. 
\end{internalclaim}
\begin{proof}
We consider a part of the graph of the function $\BasicEdgepathSystem_\dec$ in the region $0\le u \le 1/2$.
The slope of the last segment of $\BasicEdgepathSystem_\dec$ is $R=\sum_{i=1}^{\NumTangles} (-1/r_i)$.
The equation $\widetilde{\BasicEdgepathSystem_\dec}(u_0)=0$
is expressed by $R \cdot u_0-1=0$ in the interval $0\le u_0\le 1/2$.
In these cases, $R$ is $2$ or greater
since the final $r$-values include two $r_i$'s equal to $-1$ or include one $r_i$ equal to $-1$ and two $r_i$'s equal to $-2$.
Hence, the equation has a solution $u_0=1/R \in (0,1/2]$.
Thus, the monotonically decreasing type I edgepath system $\EdgepathSystem_{\typeI,\dec}$ exists,
and its endpoints have the common $u$-coordinates $0<u_0\le 1/2$.
By Lemma \ref{Lem:Incompressibility}(\ref{Lem:Incompressibility:MonotonicTypeI}), $\EdgepathSystem_{\typeI,\dec}$ is incompressible.
As in Lemma \ref{Lem:Twist:TypeI}(\ref{Lem:Twist:TypeI:Case2})(a1),
the twist of $\EdgepathSystem_{\typeI,\dec}$ is $\Twist(\BasicEdgepathSystem_\dec)-4 \le \Twist < \Twist(\BasicEdgepathSystem_\dec)-2$.
\end{proof}


Here, we prove the claim for the remaining Case 2-2-2-1.
In this case,
a Montesinos knot $K$ satisfies
$\NumTangles=3$,
$\BasicEdgepathSystem_\dec(0)=-1$,
$\sharp\{i\,|\,r_i=-1\}=1$,
$\sharp\{i\,|\,r_i=-2\}=1$ and
$\sharp\{i\,|\,r_i\le -3\}=1$.
Without loss of generality,
we can assume that $r_1=-1$, $r_2=-2$ and $r_3\le -3$.
Furthermore,
we can assume that
$T_1$ satisfies $-1<T_1<0$
and others satisfy $0<T_i<1$.
That is, any knot $K$ in this case can be isotoped 
into the normalized one.
We here divide further the Montesinos knots in Case 2-2-2-1 into three cases.

\begin{description}
 \item[Case a] $K$ satisfies $-1/3\le T_1<0$.
 \item[Case b] $K$ satisfies $-1/2\le T_1<-1/3$.
  \begin{description}
   \item[Case b-a] $K$ satisfies $r_3=-3$ or $-4$.
   \item[Case b-b] $K$ satisfies $r_3\le-5$.
  \end{description}
\end{description}

%
\begin{internalclaim}[Cases a, b-a] 
For a Montesinos knot in these cases,
the monotonically decreasing type I edgepath system 
$\EdgepathSystem_{\typeI,\dec}$
exists,
an essential surface is associated to $\EdgepathSystem_{\typeI,\dec}$,
and
$\EdgepathSystem_{\typeI,\dec}$ has twist $\Twist$ satisfying $\Twist(\BasicEdgepathSystem_\dec)-6\le\Twist<\Twist(\BasicEdgepathSystem_\dec)-4$.
\end{internalclaim}
\begin{proof}
$\BasicEdgepathSystem_\dec(1/2)<0$ holds
since 
$\BasicEdgepath_{\dec,1}(1/2)=-1/2$,
$\BasicEdgepath_{\dec,2}(1/2)=1/4$,
and $\BasicEdgepath_{\dec,3}(1/2)\le 1/6$ hold.
Besides, $\BasicEdgepathSystem_\dec(2/3)\ge 0$ is obtained
from $\BasicEdgepath_{\dec,1}(2/3)= -1/3$,
 $\BasicEdgepath_{\dec,2}(2/3)= 1/3$ and $\BasicEdgepath_{\dec,3}(2/3)>0$
in Case a,
from $\BasicEdgepath_{\dec,1}(2/3)\ge -1/2$, $\BasicEdgepath_{\dec,2}(2/3)= 1/3$ and $\BasicEdgepath_{\dec,3}(2/3)\ge 1/6$
in Case b-a.
Then the solution $1/2 < u_0\le 2/3$ exists for 
the equation $\widetilde{\BasicEdgepathSystem_\dec}(u)=0$. 
Hence, 
there exists
the monotonically decreasing type I edgepath system $\EdgepathSystem_{\typeI,\dec}$ 
with $1/2< u_0\le 2/3$.
By Lemma \ref{Lem:Incompressibility}(\ref{Lem:Incompressibility:MonotonicTypeI}), 
$\EdgepathSystem_{\typeI,\dec}$ is incompressible.
By Lemma \ref{Lem:Twist:TypeI}(\ref{Lem:Twist:TypeI:Case2})(a2), 
$\EdgepathSystem_{\typeI,\dec}$ has $\Twist(\BasicEdgepathSystem_\dec)-6\le\Twist<\Twist(\BasicEdgepathSystem_\dec)-4$.
\end{proof}

%
\begin{internalclaim}[Case b-b]
For a Montesinos knot in this case,
from $\BasicEdgepathSystem_\dec$,
a minimal class B type II edgepath system 
$\EdgepathSystem_{\typeII,\classB}$ 
is obtained by replacing 
a $\angleb{-1}$\,--\,$\angleb{-1/2}$ edge of $\BasicEdgepath_{\dec,1}$
by a $\angleb{0}$\,--\,$\angleb{-1/2}$ edge.
An essential surface is associated to $\EdgepathSystem_{\typeII,\classB}$,
and
$\EdgepathSystem_{\typeII,\classB}$ has twist $\Twist$ satisfying $\Twist=\Twist(\BasicEdgepathSystem_\dec)-4$.
\end{internalclaim}
\begin{proof}
%
Let $\Edgepath_{\typeII,\classB,i}$ denote the $i$-th edgepath of
$\EdgepathSystem_{\typeII,\classB}$.
Note that the edgepath $\Edgepath_{\typeII,\classB,1}$ is still minimal
since $T_1<-1/3$ holds and
$\Edgepath_{\typeII,\classB,1}$ does not include the edge $\angleb{-1/2}$\,--\,$\angleb{-1/3}$.
$\EdgepathSystem_{\typeII,\classB}$ is incompressible
by Lemma \ref{Lem:Incompressibility}(\ref{Lem:Incompressibility:TypeIorII})
since the cycle of final $r$-values is $(+1,-2,r_3)$ with $r_3\le -5$.
$\EdgepathSystem_{\typeII, \classB}$ has 
twist $\Twist=\Twist(\BasicEdgepathSystem_\dec)-4$, 
by Lemma \ref{Lem:Twist:TypeIIorIII}(2)(b).
\end{proof}

Consequently, with these six claims, 
we have shown Proposition \ref{Prop:Twist:LowerBound}.
\end{proof}


\begin{remark}
\label{Rem:SignOfMaxMinTwists}
The twists of the maximal and minimal boundary slopes
are always non-negative and non-positive respectively, as follows.
Since any basic edgepath has at least length $1$ and 
any basic edgepath system has $\NumTangles\ge 3$ edgepaths,
we have $\Twist_{\dec}\ge 6$ and $\Twist_{\inc}\le -6$.
As described in the argument above,
for the maximal or minimal boundary slope,
both $\Twist_{\max}-\Twist_\dec\ge -6$
and $\Twist_{\min}-\Twist_\inc\le 6$ hold.
Hence, $\Twist_{\min}\le 0\le \Twist_{\max}$.
\end{remark}

\begin{remark}
\label{Rem:BestPossibilityOfThm:Diam:LowerBound:ByCrossing}
Though it is not completely confirmed,
the lower bound given in Theorem \ref{Thm:Diam:LowerBound:ByCrossing} seems to be sharp in a sense.
We think about a family of Montesinos knot $K=M(-1/3,1/3,1/n)$, which correspond to Claim 5.
By calculation,
$\Diam(K)-2\,\Crossing(K)$ seems to become arbitrary close to $-6$
as $n$ goes to infinity.
However, it seems that no Montesinos knots $K$ can achieve $\Diam(K)=2\,\Crossing(K)-6$.
\end{remark}

\section{Remainder terms of edgepath systems}
\label{Sec:Remainder}

In this section, we will prove Theorem \ref{Thm:Diam:LowerBound:Euler}.
The keys are the ``remainder term'' of an edgepath system
and 
Proposition \ref{Prop:Diam:LowerBound:Main},
which is a technical proposition about the remainder term.

The comparison of the twist and the Euler characteristic;
strictly speaking, 
the ratio of the Euler characteristic and the number of sheets
for an essential surface with the maximal or minimal boundary slope, 
plays a main role in the proof of Theorem \ref{Thm:Diam:LowerBound:Euler}. 
Let $\EdgepathSystem$ be an edgepath system with twist $\Twist (\EdgepathSystem)$. 
Then, as we saw in \cite{IM1}, the ratio $-\chi/\sharp s (F)$ 
for a surface $F$ associated to $\EdgepathSystem$ is roughly 
the sum of lengths of the edgepaths in $\EdgepathSystem$. 
Remark that this ratio is determined independently 
from the choice of the surfaces corresponding to $\EdgepathSystem$. 
See \cite{IM1} for details. 
Thus we also use the notation $-\chi/\sharp s (\EdgepathSystem)$ 
for the value corresponding to a surface associated to $\EdgepathSystem$. 
On the other hand, 
by definition, $\Twist (\EdgepathSystem)$ is roughly calculated as 
twice of the sum of signed lengths of the edges in $\EdgepathSystem$. 
Thus, if $\EdgepathSystem$ is nearly monotonic,
$\Twist (\EdgepathSystem)$ is roughly twice of $-\chi/\sharp s (\EdgepathSystem)$.
In view of these, 
we define the \textit{remainder term} $\Remainder (\EdgepathSystem)$ 
of $\EdgepathSystem$ as 
$\Remainder (\EdgepathSystem)
=
|\Twist (\EdgepathSystem)|-2\,(-\chi/\sharp s  (\EdgepathSystem))
$.

In order to simplify the description, in the following, 
we use notations $\Remainder (F)$ and $\Twist (F)$ instead of 
the $\Remainder (\EdgepathSystem)$ and $\Twist (\EdgepathSystem)$, 
where $\EdgepathSystem$ is associated to the surface $F$, in case.

%
%
\subsection{Proposition for Reminder terms} 

Now let us state a technical proposition. 
The rest of this paper will be devoted to proving this proposition. 
%
\begin{proposition}
\label{Prop:Diam:LowerBound:Main}
For a Montesinos knot  
$K=M(\Tangle_1,\Tangle_2,\ldots,\Tangle_\NumTangles)$ with $\NumTangles\ge 3$, 
there exists an essential surface $F$ such that
its twist $\Twist (F)$ is maximal
and its remainder term $\Remainder (F)$ is non-negative.
This assertion also holds for an essential surface with the minimal twist.
\end{proposition}
We can deduce Theorem \ref{Thm:Diam:LowerBound:Euler} from this proposition.

\begin{proof}[Proof of Theorem \ref{Thm:Diam:LowerBound:Euler}]
For a non-two-bridge Montesinos knot $K$,
by Proposition \ref{Prop:Diam:LowerBound:Main}, 
there exist two essential surfaces for $K$
such that 
their twists $\Twist_1$ and $\Twist_2$ are maximal and minimal and 
their remainder terms $\Remainder_1$ and $\Remainder_2$ are both non-negative. 
Let $\Slope_1$ and $\Slope_2$ be the maximum and the minimum 
among non-meridional boundary slopes for $K$, respectively. 
Then, 
together with the definition of the twist, we have
\begin{eqnarray*}
|\Slope_1-\Slope_2|
&=&|\Twist_1-\Twist_2|=|\Twist_1|+|\Twist_2| 
=2\,(\frac{-\chi_1}{\sharp s_1}+\frac{-\chi_2}{\sharp s_2})+\Remainder_1+\Remainder_2 \\
&\ge& 2\,(\frac{-\chi_1}{\sharp s_1}+\frac{-\chi_2}{\sharp s_2})
.
\end{eqnarray*}
Note that $\Twist_1$ and $\Twist_2$ can be confirmed to be non-negative and non-positive as in Remark \ref{Rem:SignOfMaxMinTwists}.

Next, assume that a Montesinos knot $K$ is also a two-bridge knot.
The monotonically decreasing and monotonically increasing type III edgepaths
correspond to essential surfaces by \cite{HT}
and give the maximal and minimal slopes 
by the estimation of the twist of an edgepath in the previous section.
We can easily check that the remainder terms are $2$ for both essential surfaces.
Hence,
\[
\Diam(K)=2\,(\frac{-\chi_1}{\sharp s_1}+\frac{-\chi_2}{\sharp s_2})+4.
\]
\end{proof}

%
%
\subsection{Calculation of the remainder term} 

Here we include a lemma about 
calculations of the remainder terms of type I, II and III edgepath systems. 

We fix an edgepath system $\EdgepathSystem$ 
with the twist $\Twist (\EdgepathSystem)$.
We recall one more definition used in \cite[Subsection 5.2]{IM1}. 
We collect all non-$\infty$-edges $\Edge_{i,j}$ of all non-constant edgepaths $\Edgepath_i$
in $\EdgepathSystem$,
divide them into two classes 
according to the sign $\sigma(\Edge_{i,j})$,
and then sum up the lengths of edges for each class. 
With the total lengths $l_{+}$ and $l_{-}$,
let $\Cancel(\EdgepathSystem)$ denote $\min (l_{+}, l_{-})$, and 
call it the \textit{cancel} of the edgepath system. 
With the cancel, we have the following.

\begin{lemma}
\label{Lem:RemainderTerm}
The remainder term $\Remainder$ 
of an edgepath system $\EdgepathSystem$ is calculated
with the cancel $\Cancel$
as follows.
\begin{enumerate}

 \item \label{Lem:RemainderTerm:TypeI}
  Assume that 
  $\EdgepathSystem$ is a type I edgepath system.
  Let $\ConstantEdgepaths$ and $\NumTangles_{\mathrm{const}}$
  denote the constant edgepaths in $\EdgepathSystem$
  and the number of the constant edgepaths,
  where each constant edgepath $\Edgepath_i$ is a point
  on a horizontal edge $\angleb{p_i/q_i}$\,--\,$\circleb{p_i/q_i}$.
  Then, the remainder term is calculated by
  \[ 
    \Remainder
    = 
      -4 \Cancel
      +2\,(\NumTangles-\NumTangles_{\mathrm{const}})
      -\left(
       \NumTangles-2-\sum_{\Edgepath_i \in \ConstantEdgepaths}\frac{1}{q_i}
       \right)
       \frac{2}{1-u}
       .
  \]

 \item \label{Lem:RemainderTerm:TypeII}
  For a type II edgepath system $\EdgepathSystem$,
  we have $\Remainder=4-4\Cancel$.
  Thus, $\Remainder\ge 0$ holds if $\Cancel\le 1$.

 \item \label{Lem:RemainderTerm:TypeIII}
  For a ``non-augmented'' type III edgepath system $\EdgepathSystem$
  with complete $\infty$-edges, 
  we have $\Remainder=-4\Cancel$.
  Thus, $\Remainder\ge 0$ holds if $\Cancel=0$,
  that is, the type III edgepath system is monotonic.
\end{enumerate}
\end{lemma}

\begin{proof}
For non-$\infty$-edges of an edgepath system,
let $A$ be the sum of length of the edges,
$B$ the sum of signed length of the edges.
Then,
we have
$A=l_{+}+l_{-}$,
$B=l_{+}-l_{-}$,
$|B|=|l_{+}-l_{-}|=l_{+}+l_{-}-2\min(l_+,l_-)=A-2\Cancel$,
$\Twist=-2B$ and
$|\Twist|=2|B|=2\,(A-2\,\Cancel)=2\,A-4\,\Cancel$.

(\ref{Lem:RemainderTerm:TypeI})
It is just a calculation of $\Remainder = |\Twist|-2\,(-\chi/\sharp s)$ with 
the following formula in \cite{IM1}.
\begin{eqnarray*}
\frac{-\chi}{\sharp s}
&=&
 A
+\NumTangles_{\mathrm{const}}-\NumTangles
+\left(
\NumTangles-2-\sum_{\Edgepath_i \in \ConstantEdgepaths}\frac{1}{q_i}
\right)
\frac{1}{1-u}
.
\nonumber
\end{eqnarray*}

(\ref{Lem:RemainderTerm:TypeII})
By the formula $-\chi/\sharp s = A-2$ given in \cite{IM1},
we have $\Remainder = |\Twist|-2\,(-\chi/\sharp s) = 4-4\Cancel$.
%

(\ref{Lem:RemainderTerm:TypeIII})
By $-\chi/\sharp s = A$ given in \cite{IM1},
we have $\Remainder = |\Twist|-2\,(-\chi/\sharp s)=-4\Cancel$.
\end{proof}
Note that the term ``augmented'' is mentioned in the last subsection.

%
%
\subsection{Proof of Proposition \ref{Prop:Diam:LowerBound:Main}} 
This subsection is devoted to proving Proposition \ref{Prop:Diam:LowerBound:Main}. 

Our strategy to prove the proposition is as follows: 
In the proof of Proposition \ref{Prop:Twist:LowerBound},
we divided the set of Montesinos knots into some classes, and 
found for each class an edgepath system whose twist gives a lower bound on the maximal twist. 
As in that proof, we will divide the set of Montesinos knots, and 
for a Montesinos knot in each class, 
collect edgepath systems with the twist equal to or greater than that lower bound, 
and then prove that their reminder terms are all non-negative.

We first claim that 
the arguments for some special classes of edgepath systems can be omitted. 

\begin{lemma}
\label{Lem:SpecialEdgepathSystem1}
The edgepath system $\EdgepathSystem$ 
for a Montesinos knot which is either 
\begin{enumerate}
 \item[(a)] augmented type III edgepath systems,
 \item[(b)] type III edgepath systems with partial $\infty$-edges, or 
 \item[(c)] type II edgepath systems with redundant vertical edges 
\end{enumerate}
satisfies either of the following conditions (i), (ii), (iii) or (iv).

Assume that $F_1$ denotes an essential surface corresponding to $\EdgepathSystem$. 
Then, 
\begin{enumerate}
 \item[(i)]
  $\Twist(F_1)=\Twist(F_2)$ and
  $(-\chi/\sharp s)(F_1)\ge(-\chi/\sharp s)(F_2)$ hold
  for some essential surface $F_2$.
 \item[(ii)]
  $\Twist(F_1)<\Twist(F_2)$ holds for some essential surface $F_2$.
 \item[(iii)]
  $\Remainder(F_1)\ge 0$ holds.
 \item[(iv)]
  $\Twist(F_1)=\Twist(F_2)$ and $\Remainder(F_2)\ge 0$ hold
  for some essential surface $F_2$.
\end{enumerate}
\end{lemma}

By virtue of this lemma, in the following proof of the proposition, 
we can ignore edgepath systems satisfying; 
the condition (i), for it is sufficient to check if $\Remainder(F_2)\ge 0$; 
the condition (ii), for $\Twist(F_1)$ cannot be the maximum; 
the conditions (iii) and (iv), 
for even if $F_1$ gives the maximal slope, it gives $\Remainder\ge 0$. 
  
We prepare another lemma as follows. 

\begin{lemma}
\label{Lem:SpecialEdgepathSystem2}
Let $K$ be a Montesinos knot such that 
the monotonically decreasing basic edgepath system $\BasicEdgepathSystem_\dec$ 
satisfies $\BasicEdgepathSystem_\dec(0)=-1$ or $0$. 
Assume that $F_1$ denotes an essential surface with twist $\Twist(F_1)$ 
associated to a class B or class C type III edgepath system for $K$. 
Then 
there exists an essential surface $F_2$ with twist $\Twist(F_2)$ 
associated to the monotonically decreasing type III edgepath system 
such that $\Twist(F_1)<\Twist(F_2)$ holds. 
\end{lemma}

Since their proofs are rather technical, 
we give them in the next subsection separately, 
in order to make the arguments simpler.

\begin{proof} [Proof of Proposition \ref{Prop:Diam:LowerBound:Main}]
Following the strategy as we stated above, let us start to prove of the proposition.

\setcounter{internalclaimcounter}{0}


\begin{internalclaim}[Case 1]
For a Montesinos knot satisfying $\BasicEdgepathSystem_\dec(0) \ge 0$,
there exists an essential surface 
associated to a monotonically decreasing type II edgepath system 
such that its twist is maximal
and its remainder term is non-negative.
\end{internalclaim}
\begin{proof}
As in the proof of Proposition \ref{Prop:Twist:LowerBound}, 
there exists an incompressible, 
a monotonically decreasing type II edgepath system $\EdgepathSystem_{\typeII,\dec}$. 
Its twist $\Twist(\EdgepathSystem_{\typeII,\dec})$ is maximal 
by Lemmas \ref{Lem:Twist:TypeI}(\ref{Lem:Twist:TypeI:General}) 
and \ref{Lem:Twist:TypeIIorIII}.
Its remainder term $\Remainder (\EdgepathSystem_{\typeII,\dec})$ is 
$4$ by Lemma \ref{Lem:RemainderTerm}(\ref{Lem:RemainderTerm:TypeII}).
\end{proof}

%
\begin{internalclaim}[Case 3]
For a Montesinos knot satisfying $\BasicEdgepathSystem_\dec(0) \le -2$,
there exists an essential surface 
associated to the monotonically decreasing type III edgepath system 
such that its twist is maximal
and its remainder term is non-negative.
\end{internalclaim}
\begin{proof}
As in the proof of Proposition \ref{Prop:Twist:LowerBound},
the monotonically decreasing type III edgepath system 
$\EdgepathSystem_{\typeIII,\dec}$ is incompressible. 
Its twist $\Twist(\EdgepathSystem_{\typeIII,\dec})$ is maximal 
by Lemmas \ref{Lem:Twist:TypeI}(\ref{Lem:Twist:TypeI:General}) 
and \ref{Lem:Twist:TypeIIorIII}.
Its remainder term $\Remainder (\EdgepathSystem_{\typeIII,\dec})$ is 
$0$ by Lemma \ref{Lem:RemainderTerm}(\ref{Lem:RemainderTerm:TypeIII}).
\end{proof}

%
\begin{internalclaim}[Cases 2-1, 2-2-3, 2-3-2]
%
For a Montesinos knot in these cases, 
there exists an essential surface associated to 
a class A type II edgepath system $\EdgepathSystem_{\typeII,\classA}$ or 
the monotonically decreasing type III edgepath system 
such that its twist is maximal
and its remainder term is non-negative.
\end{internalclaim}
\begin{proof}
As in the proof of Proposition \ref{Prop:Twist:LowerBound},
there exists an incompressible, 
class A type II edgepath system $\EdgepathSystem_{\typeII,\classA}$. 
Its twist $\Twist (\EdgepathSystem_{\typeII,\classA})$ 
is maximal among those of type I or type II edgepath systems 
by Lemmas \ref{Lem:Twist:TypeI}(\ref{Lem:Twist:TypeI:General}) and 
\ref{Lem:Twist:TypeIIorIII}(\ref{Lem:Twist:TypeIIorIII:TypeII}). 
Its remainder term $\Remainder (\EdgepathSystem_{\typeII,\classA})$ is $0$ 
by Lemma \ref{Lem:RemainderTerm}(\ref{Lem:RemainderTerm:TypeII}), 
since its cancel $\Cancel (\EdgepathSystem_{\typeII,\classA})$ is equal to $1$. 

Only the monotonically decreasing type III edgepath system 
$\EdgepathSystem_{\typeIII,\dec}$
can have the twist greater than $\Twist(\EdgepathSystem_{\typeII,\classA})$ 
by Lemma \ref{Lem:Twist:TypeIIorIII}(\ref{Lem:Twist:TypeIIorIII:TypeIII}). 
Its remainder term $\Remainder (\EdgepathSystem_{\typeIII,\dec})$ is equal to $0$ 
by Lemma \ref{Lem:RemainderTerm}(\ref{Lem:RemainderTerm:TypeIII}).
\end{proof}

%
\begin{internalclaim}[Cases 2-2-1, 2-2-2-2, 2-3-1]
%
For a Montesinos knot in these cases,
there exists an essential surface associated to 
a monotonically decreasing type I edgepath system or 
a class A type II or III edgepath system 
such that its twist is maximal
and its remainder term is non-negative.
\end{internalclaim}
\begin{proof}
As in the proof of Proposition \ref{Prop:Twist:LowerBound}, 
there exists an incompressible 
monotonically decreasing type I edgepath system $\EdgepathSystem_{\typeI,\dec}$. 
Its twist $\Twist (\EdgepathSystem_{\typeI,\dec})$ 
satisfies 
$\Twist(\BasicEdgepathSystem_\dec) -4
\le \Twist (\EdgepathSystem_{\typeI,\dec})
<\Twist(\BasicEdgepathSystem_\dec)-2$ and 
is maximal among those of type I edgepath systems 
by Lemma \ref{Lem:Twist:TypeI}(\ref{Lem:Twist:TypeI:Case2})(a1). 
%
%
%
%
\begin{subclaim}
Its remainder term $\Remainder (\EdgepathSystem_{\typeI,\dec})$ 
is non-negative. 
\end{subclaim}
\begin{proof}
The piecewise linear equation (\ref{Eq:EquationGluingConsistency})
has the form $R\cdot u-1=0$ 
for $0<u<1/2$
where $R=\sum_{i=1}^{\NumTangles} 1/(-r_i)$.
In this case,
$\sharp\{i\,|\,r_i=-1\}$ is equal to or greater than $\sharp\{i\,|\,r_i\le -3\}$.
Other $r$-values are all $-2$.
Thus, the mean value of $\parestart -1/r_1, -1/r_2, \ldots, -1/r_\NumTangles\pareend$ is $1/2$ or greater.
Hence, we have $R\ge \NumTangles/2$ and a solution $u_0=1/R\le 2/\NumTangles$.
Note that $\EdgepathSystem_{\typeI,\dec}$ does not include any constant edgepath since $u_0<1/2$ holds as shown in the proof of Proposition \ref{Prop:Twist:LowerBound}.
Eventually, by Lemma \ref{Lem:RemainderTerm}(\ref{Lem:RemainderTerm:TypeI}),
$\Remainder=2\NumTangles-(\NumTangles-2)\cdot 2/(1-u)
\ge 2\NumTangles-(\NumTangles-2)\cdot 2\NumTangles/(\NumTangles-2) = 0$.
\end{proof}

Only class A type II edgepath systems can have 
the twist greater than $\Twist(\EdgepathSystem_{\typeI,\dec})$, 
which is $\Twist(\BasicEdgepathSystem_\dec)-2$, 
among all type II edgepath systems 
by Lemma \ref{Lem:Twist:TypeIIorIII}(\ref{Lem:Twist:TypeIIorIII:TypeII}). 
Its remainder term is equal to $0$ 
by Lemma \ref{Lem:RemainderTerm}(\ref{Lem:RemainderTerm:TypeII}). 

Only class A type III edgepath system can have 
the twist greater than $\Twist(\EdgepathSystem_{\typeI,\dec})$, 
which is $\Twist(\BasicEdgepathSystem_\dec)$, 
among all type III edgepath systems 
by Lemma \ref{Lem:Twist:TypeIIorIII}(\ref{Lem:Twist:TypeIIorIII:TypeIII}). 
Its remainder term is equal to $0$ 
by Lemma \ref{Lem:RemainderTerm}(\ref{Lem:RemainderTerm:TypeIII}). 
\end{proof}

%
\begin{internalclaim}[Cases a, b-a] 
%
For a Montesinos knot in these cases,
there exists an essential surface associated to 
a monotonically decreasing type I edgepath system, 
a class A or class B type II edgepath system, or 
the monotonically decreasing type III edgepath system 
such that its twist is maximal
and its remainder term is non-negative.
\end{internalclaim}
\begin{proof}
As in the proof of Proposition \ref{Prop:Twist:LowerBound}, 
there exists an incompressible 
monotonically decreasing type I edgepath system $\EdgepathSystem_{\typeI,\dec}$. 
Its twist $\Twist (\EdgepathSystem_{\typeI,\dec})$ satisfies 
$\Twist(\BasicEdgepathSystem_\dec) -6
\le \Twist (\EdgepathSystem_{\typeI,\dec})
<\Twist(\BasicEdgepathSystem_\dec)-4$ and 
is maximal among those of type I edgepath systems 
by Lemma \ref{Lem:Twist:TypeI}(\ref{Lem:Twist:TypeI:Case2})(a2).
%
%
%
\begin{subclaim}
Its remainder term $\Remainder (\EdgepathSystem_{\typeI,\dec})$ is non-negative. 
\end{subclaim}
\begin{proof}
Recall that in this case
$\NumTangles=3$ holds and the final $r$-values are $(-1,-2,r_3)$ with $r_3\le -3$.
Since the solution $u_0$ satisfies $u_0\le 2/3$,
at most one constant edgepath system exists in $\EdgepathSystem_{\typeI,\dec}$,
and if exists, the constant edgepath is on the edge $\angleb{-1/2}$\,--\,$\circleb{-1/2}$.
By Lemma \ref{Lem:RemainderTerm}(\ref{Lem:RemainderTerm:TypeI}),
if $\NumTangles_\const=1$, then $\Remainder=4-1/(1-u_0)\ge 1$.
Otherwise, $\Remainder =6-2/(1-u_0)\ge 0$.
\end{proof}

Only class A or class B type II edgepath systems 
can have the twist greater than $\Twist(\EdgepathSystem_{\typeI,\dec})$ 
among type II edgepath systems 
by Lemma \ref{Lem:Twist:TypeIIorIII}(\ref{Lem:Twist:TypeIIorIII:TypeII}). 
Their remainder terms are $0$ 
by Lemma \ref{Lem:RemainderTerm}(\ref{Lem:RemainderTerm:TypeII}), 
since their cancel $\Cancel (\EdgepathSystem_{\typeII,\classA})$ are equal to $1$. 

Only the monotonically decreasing, class B, or class C 
type III edgepath systems can have 
the twist greater than $\Twist(\EdgepathSystem_{\typeI,\dec})$ 
among all type III edgepath systems 
by Lemma \ref{Lem:Twist:TypeIIorIII}(\ref{Lem:Twist:TypeIIorIII:TypeIII}). 
However, in this case, class B or class C type III edgepath systems 
cannot have the maximal twist by Lemma \ref{Lem:SpecialEdgepathSystem2}. 
%
The remainder term of 
the monotonically decreasing type III edgepath systems 
is equal to $0$ 
by Lemma \ref{Lem:RemainderTerm}(\ref{Lem:RemainderTerm:TypeIII}). 
\end{proof}

%
\begin{internalclaim}[Case b-b]
%
For a Montesinos knot in this case,
there exists an essential surface associated to 
a class A or class B type II edgepath system or 
the monotonically decreasing type III edgepath system 
such that its twist is maximal
and its remainder term is non-negative.
\end{internalclaim}
\begin{proof}
As in the proof of Proposition \ref{Prop:Twist:LowerBound},
there exists an incompressible 
class B type II edgepath system $\EdgepathSystem_{\typeII,\classB}$. 
Its twist $\Twist (\EdgepathSystem_{\typeII,\classB})$ is equal to 
$\Twist(\BasicEdgepathSystem_\dec)-4$. 
This is maximal among those of type I edgepath systems 
by Lemma \ref{Lem:Twist:TypeI}(\ref{Lem:Twist:TypeI:Case2})(b), 
for $\BasicEdgepath_\dec(1/2)<0$ is obtained from
$\BasicEdgepath_{\dec,1}(1/2)=-1/2$,
$\BasicEdgepath_{\dec,2}(1/2)=1/4$ and
$\BasicEdgepath_{\dec,3}(1/2)=1/(-2r_3)\le 1/10$. 
Its remainder term is $0$ 
by Lemma \ref{Lem:RemainderTerm}(\ref{Lem:RemainderTerm:TypeII}), 
since its cancel $\Cancel (\EdgepathSystem_{\typeII,\classB})$ is equal to $1$. 

Only a class A type II edgepath system 
can have the twist, which is $\Twist(\BasicEdgepathSystem_\dec)-2$, 
greater than $\Twist (\EdgepathSystem_{\typeII,\classB})$ 
among all type II edgepath systems 
by Lemma \ref{Lem:Twist:TypeIIorIII}(\ref{Lem:Twist:TypeIIorIII:TypeII}). 
Its remainder term is equal to $0$ 
by Lemma \ref{Lem:RemainderTerm}(\ref{Lem:RemainderTerm:TypeII}). 

Only the monotonically decreasing type III edgepath system 
can have the twist, which is $\Twist(\BasicEdgepathSystem_\dec)$, 
greater than $\Twist(\EdgepathSystem_{\typeII,\classB})$, 
among all type III edgepath systems 
by Lemma \ref{Lem:Twist:TypeIIorIII}(\ref{Lem:Twist:TypeIIorIII:TypeIII}).
%
Its remainder term is equal to $0$ 
by Lemma \ref{Lem:RemainderTerm}(\ref{Lem:RemainderTerm:TypeIII}). 
\end{proof}

These six claims are sufficient to prove the proposition.
\end{proof}

\begin{remark}
\label{Rem:Optimal}
We can easily confirm that the lower bound (\ref{Eq:Diam:LowerBound:Main}) of the diameter is best possible.
For example, 
assume that,
a Montesinos knot $K$ has $4$ tangles,
and its monotonically increasing and decreasing basic edgepath systems satisfy
$\BasicEdgepathSystem_\inc(0)=+2$ and $\BasicEdgepathSystem_\dec(0)=-2$.
As in Claim 2 in the proof of Proposition \ref{Prop:Diam:LowerBound:Main}, 
the monotonically decreasing type III edgepath system
and the monotonically increasing type III edgepath system
give the maximal and the minimal twists,
and both have $\Remainder=0$.
Thus, the diameter for $K$
satisfies the equality in the inequality (\ref{Eq:Diam:LowerBound:Main}).
\end{remark}

%
%
\subsection{Special edgepath systems} 

Though the argument in this subsection is necessary,
it is technical and a kind of supplement.
%
The precise definitions of special edgepaths 
which we have treated separately in the previous subsection 
are as follows:
\begin{itemize}
\item
There is an edge $\circleb{1/0}$\,--\,$\angleb{1/0}$
called the \textit{augmented edge}.
For a type III edgepath system,
in some cases,
the augmented edge can be attached to some of edgepaths.
The edgepath system thus obtained is called
an \textit{augmented type III edgepath system}.
\item
For a basic edgepath system $\BasicEdgepathSystem$
satisfying $\BasicEdgepathSystem(0)=0$,
we can attach partial $\infty$-edges with the common $u$-coordinate of ending points
to edgepaths in $\BasicEdgepathSystem$
on making a type III edgepath system.
The edgepath system thus obtained is called 
a \textit{type III edgepath system with partial $\infty$-edges}.
\item
A type II edgepath system is called a \textit{type II edgepath system with redundant vertical edges}
if it includes both upward vertical edges and downward vertical edges.
\end{itemize}

\begin{proof}[Proof of Lemma \ref{Lem:SpecialEdgepathSystem1}]

We prove the following claims one by one.
\setcounter{internalclaimcounter}{0}

\begin{internalclaim}
An edgepath system in (a) satisfies the condition (i).
\end{internalclaim}
\begin{proof}
If there is an essential surface $F_1$ 
associated to an augmented type III edgepath system, 
as mentioned in \cite{HO},
some essential surface $F_2$ is associated to 
the corresponding non-augmented type III edgepath system.
Thus $F_1$ and $F_2$ satisfy (i).
\end{proof}

\begin{internalclaim}
An edgepath system in (b) satisfies either of the condition (ii), (iii) or (iv).
\end{internalclaim}
\begin{proof}
Assume that a surface $F_1$ corresponding to a type III edgepath system with partial $\infty$-edges is essential.
Let $\BasicEdgepathSystem$ be the basic edgepath system corresponding to $F_1$.
If $\BasicEdgepathSystem$ is monotonically decreasing, then $\Remainder(F_1)\ge 0$ holds.
That is, $F_1$ satisfies (iii).
If $\BasicEdgepathSystem$ is of class B or class C,
since $\BasicEdgepathSystem_\dec(0)=-1$ or $0$ holds,
by Lemma \ref{Lem:SpecialEdgepathSystem2},
an essential surface $F_2$ is associated to the monotonically decreasing type III edgepath system. 
Then, $F_1$ and $F_2$ satisfy (ii).
If $\Twist(F_1)=\Twist(\BasicEdgepathSystem_\dec)-6$, 
then $F_1$ may give the maximal twist in Claim 5.
Though, then, at the same time, 
the monotonically decreasing type I edgepath system 
has the twist equal to $\Twist(\BasicEdgepathSystem_\dec)-6$. 
Even in this case, it will be proved that
$\Remainder\ge 0$ holds for a type I edgepath system with remainder term $\Remainder$.
Thus $F_1$ satisfies (iv).
If $\BasicEdgepathSystem$ is any other basic edgepath system,
since $\Twist(F_1)\le \Twist(\BasicEdgepathSystem_\dec)-8$ holds,
$F_1$ cannot give the maximal twist and satisfies (ii).
\end{proof}

%
%
\begin{internalclaim}
An edgepath system in (c) satisfies the condition (i).
\end{internalclaim}
\begin{proof}
Assume that a type II edgepath system $\EdgepathSystem_{\typeII,1}$ has redundant vertical edges
and that $\EdgepathSystem_{\typeII,1}$ is constructed from 
some type II edgepath system $\EdgepathSystem_{\typeII,2}$ without redundant vertical edges
by adding upward and downward vertical edges.
Assume further that an essential surface $F_1$ is associated to $\EdgepathSystem_{\typeII,1}$.
By combining the latter half of the proof of Proposition 2.9 in \cite{HO} and Exercise just after the proof,
we can see that 
if all surfaces corresponding to $\EdgepathSystem_{\typeII,2}$ are not-$\pi_1$-injective,
then all surfaces corresponding to $\EdgepathSystem_{\typeII,1}$ are also not-$\pi_1$-injective.
Since we have a non-minimal part in a deformed edgepath system in the argument, 
not-$\pi_1$-injective surfaces obtained are moreover compressible in fact.
Eventually, if no essential surface is associated to $\EdgepathSystem_{\typeII,2}$,
then no essential surface is associated to $\EdgepathSystem_{\typeII,1}$.
Now, by assumption,
$F_1$ associated to $\EdgepathSystem_{\typeII,1}$
is essential.
Hence,
there exists an essential surface $F_2$
associated to $\EdgepathSystem_{\typeII,2}$.
$\Twist(\EdgepathSystem_{\typeII,1})=\Twist(\EdgepathSystem_{\typeII,2})$
holds
since contributions of the redundant vertical edges cancel out each other.
Moreover, 
$(-\chi/\sharp s)(F_1)\ge (-\chi/\sharp s)(F_2)$ holds
since $F_2$ is simpler by the effect of redundant vertical edges.
\end{proof}

These complete the proof of the lemma. 
\end{proof}

\begin{proof}[Proof of Lemma \ref{Lem:SpecialEdgepathSystem2}]
             Let $\EdgepathSystem_{\typeIII,1}$ and $\EdgepathSystem_{\typeIII,2}$
             be a class B or class C type III edgepath system and
             the monotonically decreasing type III edgepath system
             respectively. 
             Assume that $\EdgepathSystem_{\typeIII,2}(0)=-1$ or $0$.
             Assume that an essential surface $F_1$ is associated to $\EdgepathSystem_{\typeIII,1}$.

The compressibility of a type III edgepath system
is determined by use of completely reversibility of edgepaths.
As in Lemma \ref{Lem:Incompressibility}(\ref{Lem:Incompressibility:TypeIII}),
a type III edgepath system $\EdgepathSystem$ is compressible
if and only if 
$|\EdgepathSystem(0)|\le 1$ holds
and at least $\NumTangles-2$ edgepaths are completely reversible.
Now, 
since $\EdgepathSystem_{\typeIII,2}(0)=-1$ or $0$,
the first condition holds for
both $\EdgepathSystem_{\typeIII,1}$ and $\EdgepathSystem_{\typeIII,2}$.

Suppose that
$\EdgepathSystem_{\typeIII,2}=\parestart $%
$\Edgepath_{\typeIII,2,1}$,
$\Edgepath_{\typeIII,2,2}$, $\ldots,$
$\Edgepath_{\typeIII,2,\NumTangles}$
$\pareend$ 
is compressible
and that
$\EdgepathSystem_{\typeIII,1}=\parestart $%
$\Edgepath_{\typeIII,1,1}$,
$\Edgepath_{\typeIII,1,2}$, $\ldots,$
$\Edgepath_{\typeIII,1,\NumTangles}$
$\pareend$
is not compressible.
Then, 
at least $\NumTangles-2$ edgepaths of $\EdgepathSystem_{\typeIII,2}$ 
and
at most $\NumTangles-3$ edgepaths of $\EdgepathSystem_{\typeIII,1}$
are completely reversible.
The difference appears in exactly one pair of edgepaths 
$\Edgepath_{\typeIII,2,i}$ and $\Edgepath_{\typeIII,1,i}$.
Assume that $\Edgepath_{\typeIII,2,i}$ is completely reversible
and $\Edgepath_{\typeIII,1,i}$ is not so.
Since $\Edgepath_{\typeIII,2,i}$ is
a monotonically-decreasing completely-reversible edgepath,
the edgepath is limited to an edgepath of the form
$\angleb{1/0}$\,--\,$\angleb{z}$\,--\,$\angleb{z+1/2}$%
\,--\,$\ldots$\,--\,$\angleb{z+(p-2)/(p-1)}$\,--\,$\angleb{z+(p-1)/p}$
for some integer $z$ and $p\ge 2$.
Only the other possible minimal type III edgepath system 
for the $(z+(p-1)/p)$-tangle
is
$\Edgepath=\angleb{1/0}$\,--\,$\angleb{z+1}$\,--\,$\angleb{z+(p-1)/p}$.
This edgepath $\Edgepath$ cannot be a class C type III edgepath.
For this edgepath $\Edgepath$ to be a class B type III edgepath,
$p$ must be $2$.
Though, 
$\angleb{1/0}$\,--\,$\angleb{z+1}$\,--\,$\angleb{z+1/2}$
is also completely reversible.
Thus,
if $\EdgepathSystem_{\typeIII,2}$ is compressible,
then $\EdgepathSystem_{\typeIII,1}$ is also compressible.
Similarly to the argument for (c), 
a not-$\pi_1$-injective surface obtained by Proposition 2.5 in \cite{HO}
is moreover compressible.
Eventually
if no essential surface is associated to $\EdgepathSystem_{\typeIII,2}$,
then no essential surface is associated to $\EdgepathSystem_{\typeIII,1}$.
%
Now, by assumption,
$F_1$ associated to $\EdgepathSystem_{\typeIII,1}$
is essential.
Hence,
there exists an essential surface $F_2$
associated to $\EdgepathSystem_{\typeIII,2}$.
Obviously, we have $\Twist(\EdgepathSystem_{\typeIII,1})<\Twist(\EdgepathSystem_{\typeIII,2})$.

\end{proof}


\end{document}